%% file: article_Contact_AMR_HPC.tex
\documentclass[nopreprintline,12pt]{elsarticle}

\usepackage{amssymb,amsmath}
\usepackage{hyperref}

\usepackage{caption}
\usepackage{subfig}
\usepackage{algorithm}
\usepackage{algcompatible}
\usepackage{cancel}
\usepackage{array,multirow,makecell}
\usepackage[svgnames,dvipsnames]{xcolor}
\setcellgapes{1pt}
\makegapedcells
\newcolumntype{R}[1]{>{\raggedleft\arraybackslash }b{#1}}
\newcolumntype{L}[1]{>{\raggedright\arraybackslash }b{#1}}
\newcolumntype{C}[1]{>{\centering\arraybackslash }b{#1}}

\graphicspath{{./figures}}

\usepackage{soul}
\soulregister\cite7
\soulregister\eqref7
\soulregister\ref7

\algrenewcommand\algorithmicindent{0.6em}

\begin{document}

\begin{frontmatter}



  \title{Parallel simulation and adaptive mesh refinement for 3D
    elastostatic contact mechanics problems between deformable bodies}


\author{EPALLE Alexandre$^{1,2}$, RAMIÈRE Isabelle$^1$, LATU Guillaume$^1$, LEBON Frédéric$^2$}

\address{%
$^1$ CEA, DES, IRESNE, DEC, SESC ; Cadarache, F-13108 Saint-Paul-Lez Durance, France \\
$^2$ Aix-Marseille Université, CNRS, Centrale Marseille, LMA ; F-13453 Marseille cedex 13, France \\

  {\bf This preprint corresponds to the Chapter 7 of volume 61 in AAMS,
    Advances in Applied Mechanics.}
}

\begin{abstract}
  Parallel implementation of numerical adaptive mesh refinement (AMR)
  strategies for solving 3D elastostatic contact mechanics problems is
  an essential step toward complex simulations that exceed current
  performance levels. This paper introduces a scalable,
  robust, and efficient algorithm to deal with 2D and 3D elastostatics
  contact problems between deformable bodies in a finite element
  framework. The proposed solution combines a treatment of the contact
  problem by a node-to-node pairing algorithm with a penalization
  technique and a non-conforming h-adaptive refinement of
  quadrilateral/hexahedral meshes based on an estimate-mark-refine
  approach in a parallel framework.  One of the special features of
  our parallel strategy is that contact paired nodes are hosted by the
  same MPI tasks, which reduces the number of exchanges between
  processes for building the contact operator. The mesh partitioning
  introduced in this paper respects this rule and is based on an
  equidistribution of elements over processes, without any other
  constraints.  In order to preserve the domain curvature while
  hierarchical mesh refinement, super-parametric elements are
  used. This functionality enables the contact zone to be well
  detected during the AMR process, even for an initial coarse mesh and
  low-order discretization schemes.  The efficiency of our
  contact-AMR-HPC strategy is assessed on 2D and 3D Hertzian contact
  problems. Different AMR detection criteria are considered. Various
  convergence analyses are conducted. Parallel performances up to 1024
  cores are illustrated. Furthermore, memory footprint and
  preconditionners performance are analyzed.
\end{abstract}


\begin{highlights}
\item Combinaison of adaptive mesh refinement and contact in a parallel finite element framework
\item Non-conforming h-adaptive refinement on quadrilateral and hexahedral meshes
\item Super-parametric elements to deal with curved geometries during
  the AMR process
\item Mesh partitioning ensuring contact paired nodes on the same MPI tasks
\item Very satisfactory strong scalability up to 1024 cores.
\end{highlights}


\end{frontmatter}


\section{Introduction}

Efficient numerical solving of contact mechanics problems remains
challenging. These problems are locally highly non-linear and
non-smooth. They often have the particularity of
concentrating very high stress in the contact area. Thus, the
numerical simulation of this type of problems is computationally
demanding, especially if a local accurate solution is expected.

Dealing with contact between deformable solids, the finite element
method (FEM)~\cite{Ciarlet-1978} remains the most widely used numerical
method~\cite{Wriggers-2006} for solving the associated problem, although others approaches have been
more recently introduced. Let us cite for example the virtual element method
(VEM)~\cite{BeiraodaVeiga-2013,Wriggers-2023} and its extension to
contact~\cite{Wriggers-2016,Laaziri-2025} as well as the Hybrid High Order (HHO)~\cite{DiPietro-2015,DiPietro-2020} and its adaptation to (frictional) contact problems~\cite{Chouly-2020}. The scaled boundary finite element method
(SBFEM)~\cite{Song-1997}, combining the advantages of the FEM and of
the boundary element method (BEM), has also been recently applied for contact
analysis in~\cite{Hirshikesh-2021}. Finally, let us mention the work on isogeometric contact~\cite{DeLorenzis-2014} which considers contact formulations within the framework of isogeometric analysis (IGA)~\cite{Cottrell-2009}.

Several numerical strategies have been developed over the years to
deal with constraints at contact interfaces~\cite{Kikuchi-1988, Laursen-2002,
  Wriggers-2006}. When
using the FEM, the penalization method~\cite{Kikuchi-1981}, the
Lagrange multiplier method~\cite{Kikuchi-1988}, the augmented
Lagrangian method~\cite{Simo-1992}, the Nitsche
method~\cite{Wriggers-2008,Chouly-2013b} or mortar-based
formulation~\cite{Belgacem-1998} are commonly employed to design
robust and efficient algorithms. For
more details about these strategies, the reader is referred
to~\cite{Wriggers-2006}.

Furthermore, solving contact problems with accuracy requires a high
degree of mesh refinement in the contact zones. Building uniformly
refined mesh for such problems is still prohibitively computationally
expensive, especially for 3D geometries. One way to overcome this
limit is to rely on adaptive mesh refinement (AMR)
methods~\cite{Verfurth-1996}.  These techniques, which have been well
known for many years, enable the mesh to be locally adjusted to a
desired accuracy, thereby reducing the total number of degrees of
freedom compared to uniformly refined meshes.
%
Popular adaptive techniques for reducing discretization error include
h-adaptive~\cite{Demkowicz-1985,Diez-1999} (mesh-step refinement),
p-adaptive~\cite{Babuska-1982,Babuska-1992} (order of basis function
increase), hp-adaptive~\cite{Zienkiewicz-1989,Babuska-1992}
(combination of both previous) and local
multigrid~\cite{Brandt-1977-2,Bai-1987,Khadra-1996} approaches
(geometric refinement by adding local refined meshes). Adaptive mesh
refinement is generally coupled with \textit{a posteriori} error
estimators~\cite{Verfurth-1996} in order to automatically determine
the elements to be refined, and sometimes also to automatically stop
the refinement
process~\cite{Zienkiewicz-1987, Liu-2017}.\\
Many works have focused on the application of AMR techniques in the
context of contact mechanics problems, mainly in the FEM
framework. Let us cite for example~\cite{Wriggers-1998, Braess-2007,
  Gustafsson-2020, Araya-2023} whose authors proposed dedicated
\textit{a posteriori} error estimators for contact problems that have
been applied with conforming h-adaptive refinement. In particular,
\cite{Wriggers-1998} modifies the Zienkiewicz and Zhu (ZZ)
estimator~\cite{Zienkiewicz-1987} to take into account contact
contributions.
In their turn, the works of~\cite{Roda-Casanova-2018, Rademacher-2019} focus on mesh
adaptivity with non-conforming
h-refinement methods.
Adaptation of local multigrid methods for contact problems have been
studied in~\cite{Boffy-2012, Liu-2017}. Let us mention that the
authors of~\cite{Liu-2017} successfully applied the standard ZZ
estimator in each body separately as for multimaterial problems.  Most
of the above methods have only been applied on 2D geometries. It should also be
noted that some of them only concern rigid-deformable contact problems
(Signorini problems)~\cite{Braess-2007, Roda-Casanova-2018,
  Rademacher-2019, Araya-2023}.  For its part, the extension of SBFEM
to contact problems~\cite{Hirshikesh-2021} has been directly
introduced in a non-conforming h-adaptive mesh framework based on a
quadtree decomposition with a dedicated SBFEM error
indicator~\cite{Song-2018}. Concerning the VEM, the authors
in~\cite{Aldakheel-2020} provide a methodology to combine AMR and 2D
contact problems. Two types of AMR techniques are studied on 2D
contact problems: non-conforming h-adaptive meshes or conforming
locally refined Vorono{\"i} cell meshes. The mesh is locally refined
in a predefined contact zone until a geometric criterion on the
refined contact area is reached.  The extension to 3D contact problems
of the VEM with AMR based on Vorono{\"i} type elements is proposed
in~\cite{Cihan-2022}.

However, for highly singular problems such as those arising from
contact mechanics, the discretization can become very fine locally,
and the number of unknowns prohibitive even with AMR. It is therefore
also desirable to use high-performance computing (HPC) methods for
faster processing and for handling large-scale problems.  The fast solution of partial differential
equations with the FEM on adaptively refined meshes via
parallelization strategies has already been widely explored, see~\cite{Flaherty-1997, Tu-2005, Mitchell-2007, Sundar-2008,
  Burstedde-2011, Burstedde-2020} to cite but a few. A large set of these strategies lies on
refinement-tree based partitioning methods and hence non-conforming
h-adaptive approaches (with hanging nodes). In addition, it has been recently
presented in \cite{Munch-2022} an algorithm for solving hanging-node
constraints in the context of matrix-free FEM on central processing
units (CPUs) and graphics processing units (GPUs).
All these methods
and algorithms have been integrated into several scalable open-source
FEM libraries such as deal.II~\cite{Arndt-2021-art},
MFEM~\cite{Anderson-2020}, MoFEM~\cite{Kaczmarczyk-2020},
FEniCS~\cite{Alnaes-2015}, FreeFem++~\cite{Hecht-2012} or
libMesh~\cite{Kirk-2006} for example.
\\
Addressing large-scale contact problems, adapted to the multiprocessor
architecture of supercomputers, have been under study for many
years. Indeed, pioneer
works~\cite{Barboteu-2001,Dureisseix-2001,Alart-2007} introduce domain
decomposition methods and associated iterative solution algorithms
suitable to parallel computers. However, the performance of these
methods has been shown on a small number of processes, no more than
twenty or so. Penalty method applied on rigid-deformable contact
problems~\cite{Har-2003} enables to reach good performances until
sixty processors. Recent studies~\cite{Wiesner-2021,Franceschini-2022,
  Guillamet-2022,Dostal-2023-Book,Mayr-2023} propose strategies
enabling scalable calculations on more processes. For these methods,
the increase in solution accuracy (or in number of degrees of freedom)
is based on uniform mesh refinement.  Works detailed
in~\cite{Wiesner-2021,Franceschini-2022,Mayr-2023} focus on
saddle-point problems (arising with Lagrange multipliers or mortar
methods) with dedicated preconditioners. These approaches are scalable
up to a few hundreds of cores.
In~\cite{Guillamet-2022,Dostal-2023-Book}, parallel
algorithms for rigid-deformable contact, scaling up to over 1,000 cores are
described.
The proposed strategies are based on domain
decomposition methods (Dirichlet-Neumann or FETI like) for fulfilling contact conditions.


Combining the parallelization of AMR methods with contact mechanics
remains a scientific challenge. As described before, current developments
are mainly focusing either on non-parallelized adaptive methods for
contact problems, or on parallelization methods for contact problems
solved on uniformly refined meshes. To the best of our knowledge, only
the study presented in~\cite{Frohne-2016} analyzes the parallel
performance of solving a contact mechanics problem with adaptive local
mesh refinement. A conventional one-body contact problem
(rigid-deformable contact) is solved with a primal-dual active set
strategy based on Lagrange multipliers that can be seen as a
semi-smooth Newton method. This rigid-deformable contact problem
greatly simplifies the general deformable-deformable contact problem,
especially in terms of contact detection and mesh refinement. The strategy is
implemented in the deal.II FEM library. The
mesh refinement is conducted by a variant of the Kelly error
estimator~\cite{Kelly-1983}, ready for use as part of deal.II. A
hierarchical non-conforming h-adaptive refinement is then applied to linear and quadratic
hexahedral finite elements. Contact simulation results
with a scalability up to 1,024 processes are shown.


The aim of this paper is to introduce an HPC strategy for solving 3D
elastostatic contact mechanics problems with AMR, and potentially
curved contact boundaries.  As said before, several software environments
combining HPC and AMR in a FEM framework already exist. We then decided
to add contact mechanics treatment in such environment.
The MFEM
software environment was chosen for the following main reasons:
proven scalability to large numbers of cores, mature support for
hexahedrons and tetrahedrons, support for non-conforming AMR, a large
active community, and the ability to be integrated into a software stack.
Furthermore, the minimal set of tools that MFEM provides for modeling contact mechanics is valuable.  It gave us the opportunity to develop and examine parallel solutions in more complex settings, as well as evaluate the performance of combining AMR and parallelism in such configurations.
\\
The meshes being considered are based on quadrilateral/hexahedral finite
elements.  These elements are attractive for their tensor-product
structure~\cite{Belytschko-2003, Ramiere-2007,Arndt-2021-art}, their
better accuracy (compared to triangles or tetrahedra) for the same 
approximation space order~\cite{Arndt-2021-art}, and offer great
advantages for the majority of solid mechanics problems such as
contact mechanics~\cite{Liu-2017, Michel-2013}.  A so-called
hierarchical h-adaptive AMR strategy is applied~\cite{Zhu-1993,
  Rodenas-2017}. It consists in directly
subdividing the elements in areas of
interest while maintaining the same order of interpolation, contrary
to p or hp-adaptive methods that induce an increase in the order of the elements.
Mesh-step refinement techniques enable to obtain optimal convergence
rates even for problems with geometrical
singularities~\cite{Babuska-1983, Hennig-2016}
and are strongly featured in the solid mechanics community.  In
addition, hierarchical refinement generally leads to optimal refined
meshes~\cite{Zhu-1993,Rodenas-2017,Koliesnikova-2021} with respect to
the total number of degrees of freedom and allows better control of
the adaptation process through an explicit hierarchical data
structure.  When using quad/hexa elements, this AMR approach leads to
non-conforming meshes with hanging nodes~\cite{Koliesnikova-2021} for
which continuity must be ensured.  In MFEM, a parallelized algorithm
for non-conforming h-adaptive mesh refinement is already
available~\cite{Cerveny-2019} based on prolongation and restriction
operators to solve the linear system on conforming nodes only. This
ready for use mesh refinement technique is enriched by
estimate-mark-stop-refine procedures inspired
by~\cite{Koliesnikova-2021} where they proved their effectiveness in
terms of number of mesh elements versus precision. They are mainly
based on the Zienkiewicz-Zhu \textit{a posteriori}
estimator~\cite{Zienkiewicz-1987}, already implemented in the MFEM
library. This estimator is widely used in solid mechanics as it is
generic (not confined to a specific problem), cheap to compute, easy
to implement and works pretty well in
practice~\cite{https://doi.org/10.1002/nme.2980}.
In this study, only one error estimator is tested since
our main objective is the introduction of a generic parallelized
AMR strategy for contact
mechanics and not an in-depth study of error estimators.\\
The combination of AMR and contact problems highlights the problems
of hierarchically refined curved
interfaces~\cite{Wriggers-1998,Liu-2017,Aldakheel-2020} in the FEM framework.
In order to avoid
hand-made procedures
ensuring the potential contact nodes
to be located on the curved boundary during the refinement
process~\cite{Liu-2017,Koliesnikova-2021}, we propose to
use super-parametric elements in order to preserve the shape of the
initial geometry boundary, even for first-order finite element
solutions. These are elements whose degree of interpolation of the
geometry is greater than that of the basis functions of the
solution. \\
In the present paper, the contact problem is handled by a penalization algorithm. This
enables one to deal with primal variables only: it provides a way for
modelling the phenomenon under study while not increasing the number
of unknowns (as it could be using Lagrange multipliers) and not modifying the formulation in an intrusive
manner (as it could be using Nitsche's method). In addition, parallelized operators for solving constrained
systems available in MFEM can be directly employed. Let us underline
that providing scalable solvers and associated preconditioners for saddle-point
problems is still a ongoing work~\cite{Wiesner-2021,Mayr-2023,Nataf-2023}, which supports our choice of the penalty method. 
An in-depth study around the penalty coefficient in order to limit
interpenetration while controlling the computational time is
carried out in this paper. A node-to-node pairing combined to an active set strategy~\cite{Luenberger-1984} is applied for the detection of elements in
contact.
Although it may seem limiting, this node pairing strategy remains of
interest nowadays for its good modeling properties~\cite{Jin-2016,Xing-2019}.
Its combination with automatic mesh refinement makes it possible to consider its use in the case of 
small tangential slip.
Finally, the combined AMR-contact algorithm is ruled by two nested iterative
loops. The external loop concerns the AMR process while the internal
one deals with the contact solution.\\
The proposed scalable contact algorithm is based on a mesh partitioning
strategy that ensures the contact paired nodes to be on the same MPI
process.
Hence, the algorithm building the local contact stiffness matrix does not depend on communications with other processes. 
As a result, communication overheads are reduced and the code is considerably simplified.
In practice, the partitioning of the mesh is governed by an
algorithm that guarantees the load equidistribution across
processes, and that respects the distribution constraint related to the paired contact
nodes mentioned above.
This partitioning strategy offers a straightforward implementation, as
it does not consider the spatial vicinity of elements within the
solid when assigning them to processes.\\
The scalability of the
developed AMR-contact algorithm and its main functions are analyzed in
this article.  It should be noted that mesh coarsening and mesh
adaptation over time are beyond the
scope of this article and will be considered in future studies.



The rest of the paper is organized as follows. Section
\ref{sec:contact_modeling} recalls the contact mechanics problem
between deformable bodies and presents the main ingredients of the
chosen solution algorithm. Section~\ref{sec:amr_process} is devoted to the
adaptive mesh refinement strategy, focusing on special insights for
contact mechanics. The extension of the considered AMR-contact
algorithm to parallel settings is described in Section
\ref{sec:parallel_strategy}. Various numerical and performance studies carried out
on the Hertzian elastostatic contact problem (in 2D and 3D) are reported
and discussed in Section~\ref{sec:hertzian_problem}. Finally, some concluding
remarks are provided.

\section{Contact modeling} \label{sec:contact_modeling}

\subsection{Governing equations}

Two linearly elastic solids $\Omega^{s} \left( s = 1, 2 \right)
\subset \mathbb{R}^{m}$ (with $m$ the space dimension) are considered in frictionless contact on
$\Gamma_{C}^{s}$, the potential contact boundaries. A schematic
representation of the problem is given in Figure~\ref{fig:patate_contact}.

In $\Omega = \Omega^{1} \cup \Omega^{2}$ (such that $\Omega^{1} \cap
\Omega^{2} = \emptyset$), the considered elastostatic frictionless unilateral contact problem with infinitesimal strain assumption writes
\begin{equation}
\left( \mathcal{P} \right) : \begin{cases}
- \text{div} \boldsymbol{\sigma} \left( \textbf{u} \right) = \textbf{0} & \text{ in } \Omega \\ 
\boldsymbol{\varepsilon} \left( \textbf{u} \right) = \dfrac{1}{2} \left( \nabla \textbf{u} + \nabla \textbf{u}^{T} \right) & \text{ in } \Omega \\
\boldsymbol{\sigma} \left( \textbf{u} \right) = \textbf{C} : \boldsymbol{\varepsilon} \left( \textbf{u} \right) & \text{ in } \Omega \\
\text{Boundary conditions} & \text{ on } \Gamma_{BC}\\
\text{Contact condition} & \text{ on } \Gamma_{C}.
\end{cases}
\label{elastostatic_contact_problem_hpp}
\end{equation}
with $\textbf{u}$ the problem’s solution in displacement,
$\boldsymbol{\sigma}$ and $\boldsymbol{\varepsilon}$ the stress and
strain tensors respectively, and $\textbf{C}$ the fourth order elasticity tensor. \\
Boundaries
$\Gamma_{BC}$ and $\Gamma_{C}$ are defined such as $\left\{
\begin{array}{ll}
\Gamma_{BC} \cup \Gamma_{C} &= \partial \Omega \quad \text{(boundary of $\Omega$)}  \\
\Gamma_{BC} \cap \Gamma_{C} &= \emptyset \\
\Gamma_{BC} &= \Gamma_{BC}^{1} \cup \Gamma_{BC}^{2}\\
\Gamma_{C} &= \Gamma_{C}^{1} \cup \Gamma_{C}^{2}.
\end{array}
\right.$ \\
Classical boundary conditions (Dirichlet or Neumann) are applied
on $\Gamma_{BC}$ so that the Problem $\left( \mathcal{P}
\right)$ is well-posed.


\begin{center}
	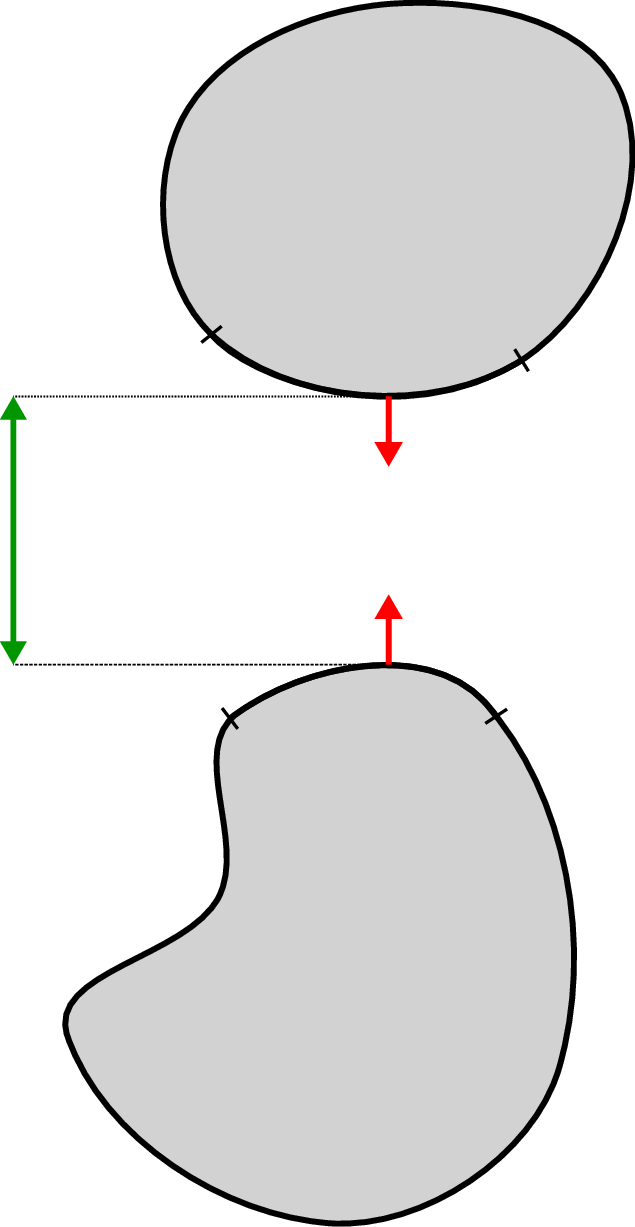
	\captionof{figure}{Sketch of the contact problem.}
	\label{fig:patate_contact}
\end{center}

In the common contact mechanics theory, Signorini contact
conditions (see~\cite{Moreau-1974, Curnier-1999}) are prescribed on
$\Gamma_{C}$, see Eq.~\eqref{eqn:signorini_conditions}.
\begin{equation} 
\begin{cases}
u_{N} - d \leq 0 \quad \text{(non penetration)}\\ 
\sigma_{N} \leq 0 \quad \text{(non adhesion)}\\ 
\sigma_{N} \left( u_{N} - d \right) = 0  \quad \text{(complementarity)}
\end{cases} \; \text{on} \; \Gamma_{C}
\label{eqn:signorini_conditions}
\end{equation}
with
$u_{N} = \textbf{u}^{1} \cdot \textbf{n}^{1} + \textbf{u}^{2} \cdot
\textbf{n}^{2}$ the normal displacement,
$\sigma_{N} = (\boldsymbol{\sigma}^{1} \cdot \textbf{n}^{1}) \cdot \textbf{n}^{1} = (\boldsymbol{\sigma}^{2} \cdot
\textbf{n}^{2}) \cdot \textbf{n}^{2} $ the normal contact stress and $d$ the initial distance
(gap) between solids. We denote by $\textbf{u}^{s}$ the displacement
field in $\Omega^{s}$: $\textbf{u}^{s} = \textbf{u}|_{\Omega^{s}} $, and by
$\boldsymbol{\sigma}^{s}$ the stress field in $\Omega^{s}$: $\boldsymbol{\sigma}^{s}= \boldsymbol{\sigma}|_{\Omega^{s}}$.
The symbol '$\cdot$' indicates the inner product,
while $\textbf{n}^{s}$ represents the external unit normal to $\Omega^{s}$ on
$\Gamma_{C}^{s}$. Moreover, the framework of the infinitesimal strain theory implies that
$\textbf{n}=\textbf{n}^{1}=-\textbf{n}^{2}$.

\begin{figure}[!h]\centering
\subfloat[][Contact law]{ \includegraphics[scale=0.5]{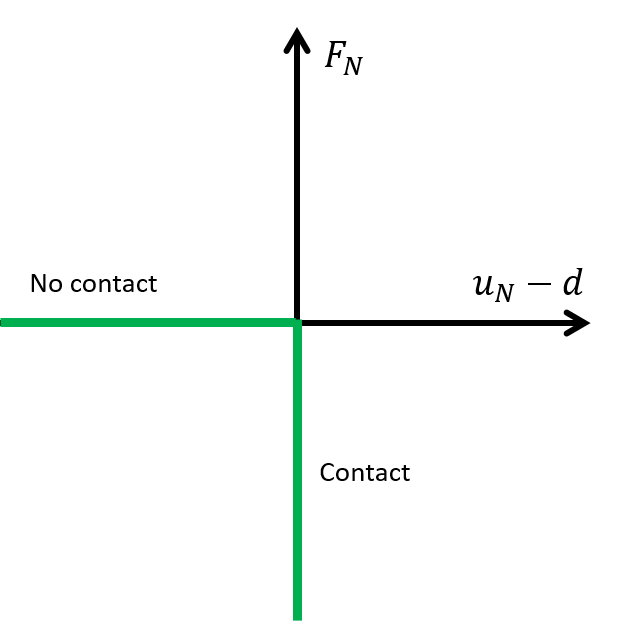} 
\label{fig:contact_figure}}
\hspace*{0.2\baselineskip} 
\subfloat[][Penalized contact law]{ \includegraphics[scale=0.5]{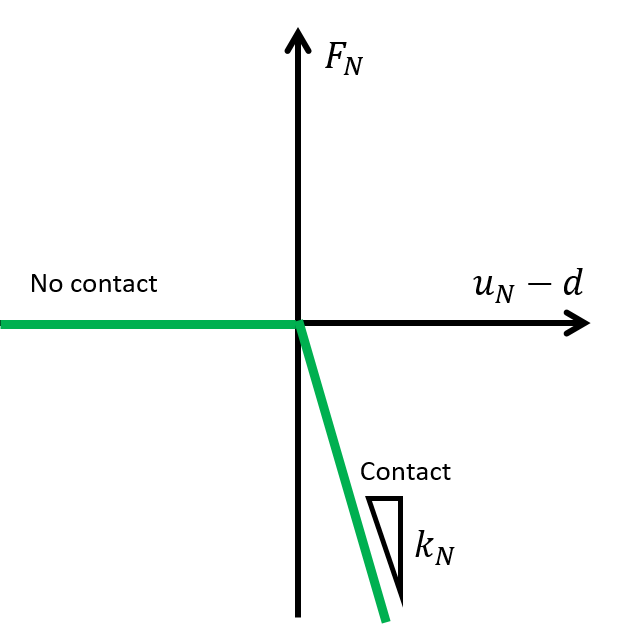}  \label{fig:contact_penalty}}
\caption{Illustration of contact and penalized contact laws.}
\label{fig:contact_laws}
\end{figure}

Conditions~\eqref{eqn:signorini_conditions}, graphically represented
in Figure~\ref{fig:contact_figure} in contact forces $F_{N}$, make the
contact problem highly non-smooth and non-linear. 

\subsection{Discrete form} \label{sec:discrete-form}

A generic contact mechanics problem in discrete form can be written as
\begin{equation}
\left( \mathcal{P}' \right) : \begin{cases}
\textbf{L}^{int} \left( \textbf{U} \right) - \textbf{L}^{ext} - \textbf{L}^{C} \left( \textbf{U} \right) &= \textbf{0} \\ 
 \textbf{B} \textbf{U} & \leq \textbf{D}
\end{cases},
\label{eqn:generic_contact_problem}
\end{equation}
where $\textbf{U}$ denotes the discrete displacement solution,
while $\textbf{L}^{int}$, $\textbf{L}^{ext}$ and
$\textbf{L}^{C}$ denote the internal, external and contact
forces operators respectively.  In
addition, $\textbf{B}$ is the pairing matrix of potential contact
nodes used to define the discrete relative displacements and $\textbf{D}$
the initial gaps
between the paired nodes.

Hence in Eq.~\eqref{eqn:generic_contact_problem}, the '$\leq$' symbol is used to represent the componentwise inequality.

\textbf{Remark:} The term $\textbf{L}^{C}(\textbf{U})$ should be understood in a general sense.
When applying the Lagrange multipliers method or one of its variants (augmented Lagrangian, mortar,...), the dependence on $\textbf{U}$ is replaced by a dependence on the additional variable associated with the method.\\

In the elastostatic case under consideration, the internal forces
reduce to
\begin{equation}
  \label{eqn:elastic-stiffness}
  \textbf{L}^{int}(\textbf{U}) = \textbf{K} \textbf{U},
\end{equation}
with $\textbf{K}$ the stiffness matrix.


The penalty method consists in rewriting the contact inequality
problem~\eqref{eqn:generic_contact_problem} as a variational equality,
much easier to solve by classical solvers. This is done by identifying the set of
active contact conditions, denoted
$\mathcal{A}$, which is defined by
Equation~\eqref{eqn:active_contact_set}.
\begin{equation}
\mathcal{A} =\{ i  \in [\![ 1;N_{C} ]\!] \; ; \; \sum_{j=1}^N \textbf{B}_{ij} \textbf{U}_{j} \geq \textbf{D}_{i}\}
\label{eqn:active_contact_set}
\end{equation}
with $N$ the number of (primal) degrees of freedom, and $N_C$ the number of
potential contacts.\\
As the matrix $\textbf{B}$ is defined only
on the degrees of freedom located the potential contact boundary $\Gamma_C$, the sum in
Eq.~\ref{eqn:active_contact_set} is practically restricted to 
the primal nodes belonging to the set $\mathcal{A}_c$
\begin{equation*}
\mathcal{A}_c =\{ j  \in [\![ 1;N ]\!] \; ; \;\exists i \in
[\![ 1;N_{C} ]\!] \; \text{such that} \; \textbf{B}_{ij} \ne 0\}.
\end{equation*}

Then, the active pairing matrix $\hat{\textbf{B}}$
is obtained by restricting the rows of $\textbf{B}$ to the set
$\mathcal{A}$, see Eq.~\eqref{eqn:pairing_matrix_active}. The
'active' gap vector $\hat{\textbf{D}}$ is defined similarly, see Eq.~\eqref{eqn:pairing_matrix_active}.
\begin{equation}
  \begin{cases}
  \hat{\textbf{B}} &= \textbf{B} \left[ \mathcal{A} \right]\\
  \hat{\textbf{D}} &= \textbf{D}\left[ \mathcal{A} \right]  
  \end{cases}
\label{eqn:pairing_matrix_active}
\end{equation}
The contact reaction force is defined in order to penalize surface
interpenetration. By the way, the contact formulation is regularized,
see Figure~\ref{fig:contact_penalty}, and the discrete penalized
contact forces $\textbf{L}^{C}(\textbf{U})$ are given by Eq.~\eqref{eqn:force-contact-penalization}.
\begin{equation} \label{eqn:force-contact-penalization}
\textbf{L}^{C}(\textbf{U}) = k_{N} \hat{\textbf{B}}^{T} \left( \hat{\textbf{D}} - \hat{\textbf{B}} \textbf{U} \right) 
\end{equation}
with $k_{N}\gg 1$ {N.m$^{-1}$} the penalty coefficient, which can be interpreted as
a spring stiffness in the contact interface, see~\cite{Wriggers-2006}. Finally the mechanical
equilibrium writes
\begin{equation}
\left( \textbf{K} + k_{N} \hat{\textbf{K}}_{C} \right) \textbf{U} = \textbf{L}^{ext} + k_{N} \hat{\textbf{B}}^{T} \hat{\textbf{D}}
\label{eqn:mechanical_equilibrium_matrix_system}
\end{equation}
with $\hat{\textbf{K}}_{C} = \hat{\textbf{B}}^{T}
\hat{\textbf{B}}$. \\

As mentioned in the introduction,
Problem~\eqref{eqn:mechanical_equilibrium_matrix_system} only involves
primal unknowns. On the other hand, with this penalty approach,
interpenetration between solids is allowed~\cite{Wriggers-2006}, since this is the
driving term of the contact force.

The \texttt{mfem::PenaltyConstrainedSolver} class of the
MFEM software fully forms the linear
system~\eqref{eqn:mechanical_equilibrium_matrix_system}
from
penalty coefficient $k_{N}$, matrices and vectors $\textbf{K}$,
$\hat{\textbf{B}}$, $\textbf{L}^{ext}$ , $\hat{\textbf{D}}$ provided as input.

\subsection{Node-to-node pairing} \label{sec:nodenode_pairing}


To solve the contact problem described by
Equation~\eqref{eqn:mechanical_equilibrium_matrix_system}, a
preliminary step which consists in pairing potential contact nodes is
required. Indeed, this step enables to build the matrix $\textbf{B}$
and the vector $\textbf{D}$, from the positions of the mesh nodes
located on the contact boundaries of the solids.
The process consists in pairing the closest potential contact nodes using a position criterion.

In computational contact mechanics, there are many ways to pair the
contact interface nodes. The most common are the node-to-node,
node-to-surface and surface-to-surface pairing algorithms. Here, the
node-to-node pairing approach is chosen for its ease of implementation
and very low calculation costs. However, it restricts the problem and
the geometry to be taken into account. The initial meshes must have
matching contact interfaces. Moreover, this pairing approach is only
valid under infinitesimal strain theory and infinitesimal sliding, so that the potential contact nodes always face each
other. 

Thus, the MFEM software is enriched with our own node-to-node pairing
algorithm which consists in finding, for each node of the first
solid, the nearest node of the second solid in the direction of
the ({\it a priori} given) contact outward normal.

\subsection{Algorithmic treatment of the penalized contact problem} \label{sec:algorithmic_treatment_contact_problem}

The proposed algorithm for solving the contact problem is ruled by an
iterative process, see
Algorithm~\ref{alg:contact_solution_process_algo}. Node-to-node
pairing
is performed before
entering the loop, whereas the set of active contact pairs of nodes
$\mathcal{A}$ is updated at the end of each iteration thanks to an
active set algorithm (see~\cite{Luenberger-1984} for example) by looking at the
interpenetration between the solids, see Definition~\eqref{eqn:active_contact_set}. The set of active contact
nodes may vary at each iteration, and the contact matrix and vector,
$\hat{\textbf{K}}_{C}$ and $\hat{\textbf{B}}^{T} \hat{\textbf{D}}$ respectively,
are reconstructed accordingly. The resulting linear system, see
Eq.~\eqref{eqn:mechanical_equilibrium_matrix_system}, is typically
solved by a preconditioned conjugate gradient (PCG) linear iterative
solver, suitable for penalized problems resulting in large sparse symmetric
positive-definite systems. The stopping criterion of the contact solution process requires that the set of active contact nodes no longer varies in iteration.
%

\begin{algorithm}[!h]
\caption{Contact solution process }
\label{alg:contact_solution_process_algo}
\begin{algorithmic}
\STATE \textbf{Input: Data of the problem, maximum number of contact
  loop iterations $l_{max}$}
\STATE \textbf{Output: Finite element solution $\textit{\textbf{U}}$}
\STATE \textbf{Initialization}
\STATE \hspace{\algorithmicindent} Process the node-to-node pairing
\STATE \hspace{\algorithmicindent} Initialize contact loop parameters: $l \leftarrow 1$, $\textit{\textbf{stop}}_{\textit{\textbf{C}}} \leftarrow false$
\STATE \hspace{\algorithmicindent} Define the initial set of active contact nodes (if any)
\WHILE{($\textit{\textbf{stop}}_{\textit{\textbf{C}}} = false$ and $l \leq l_{max}$)}
\STATE Form penalized linear system
\STATE Solve the linear system
\STATE Define the new set of active contact conditions
\IF{no evolution of the active contact set}
\STATE $\textit{\textbf{stop}}_{\textit{\textbf{C}}} \leftarrow true$ 
\ELSE
\STATE $l \leftarrow l + 1$
\ENDIF
\ENDWHILE
\end{algorithmic}
\end{algorithm}

\section{Adaptive mesh refinement process} \label{sec:amr_process}

The non-conforming h-adaptive strategy is performed using an iterative
algorithm described in Figure~\ref{fig:generic_amr_algorithm}. The
problem is first solved on a coarse (generally conforming) mesh before applying the
refinement process according to a chosen detection criterion. This process, called ESTIMATE-MARK-REFINE strategy~\cite{Verfurth-1996, Nochetto-2012, Koliesnikova-2021}, is iteratively repeated until a stopping criterion is satisfied.

\begin{center}
	\includegraphics[scale=0.15]{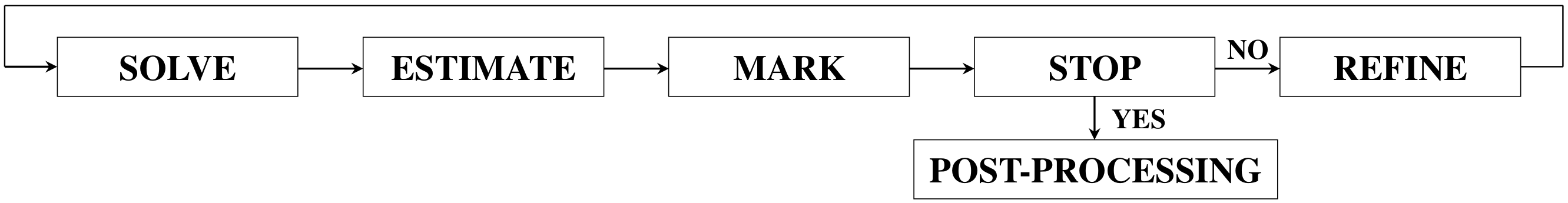} 
	\captionof{figure}{Generic adaptive mesh refinement algorithm, figure from~\cite{Koliesnikova-2021}.}
	\label{fig:generic_amr_algorithm}
\end{center}

In the following sections, the main points of each module of the AMR
algorithm are described.

\subsection{Error estimation: module ESTIMATE} \label{sec:amr_estimate}

The key ingredient for automatic mesh adaptation is the error
estimator. Since the present study does not aim at comparing the
efficiency of error estimators, the Zienkiewicz and
Zhu~\cite{Zienkiewicz-1987} {\it a posteriori} error estimator is
set as the tool for driving mesh refinement. This estimator is based
on the difference between a smoothed stress field
$\boldsymbol{\sigma}^{\star}$ and the finite element one
$\boldsymbol{\sigma}^{h}$, which is, for a displacement-based
solution, generally discontinuous between elements. For the construction of the
smoothed field, the nodal projection technique introduced
in~\cite{Zienkiewicz-1987} is performed. Assuming that $\Omega^{h}$ is
the quadrangulation of $\Omega$, each mesh element is denoted
$T_{i}, i \in [\![ 1;N_{E} ]\!]$ such that
$\overline{\Omega^{h}} = \cup_{i}^{N_{E}} \overline{T_{i}}$. Hence $N_{E}$ represents
the total number of element of the mesh. The local element-wise error
in the energy norm writes for a mesh element $T_{i}$:
\begin{equation}
\xi_{T_{i}} = \left( \int_{T_{i}} \left( \boldsymbol{\sigma}^{\star} - \boldsymbol{\sigma}^{h} \right) : \left( \boldsymbol{\varepsilon}^{\star} - \boldsymbol{\varepsilon}^{h} \right) d T_{i} \right)^{1/2} .
	\label{eqn:element-wise_error_distribution}
\end{equation}

\underline{Remark:} In practice, this integral is evaluated at
integration points where $\boldsymbol{\sigma}^{h}$ is defined.\\

As proposed in~\cite{Liu-2017}, this classical projection-based ZZ estimator~\cite{Zienkiewicz-1987} is applied in each
solid independently without any other modification. The global
absolute estimated error on $\Omega^{h}$ is then written by summing the
elementary contributions
\begin{equation}
\xi_{\Omega} = \left( \sum_{T_{i} \in \Omega^{h}} \xi_{T_{i}}^{2} \right)^{1/2} .
	\label{eqn:global_error_distribution}
\end{equation}

In the MFEM software library, the element-wise error distribution
$\{\xi_{T_{i}}\}_{T_i \in \Omega^h}$ are computed using the
\texttt{GetLocalErrors} function applied to the
\texttt{mfem::ZienkiewiczZhuEstimator} error estimator. The scalar value
$\xi_{\Omega}$ is then obtained from the element-wise field $\{\xi_{T_{i}}\}_{T_i \in \Omega^h}$ thanks to the \texttt{GetNorm} function.    

\subsection{Marking stage: module MARK} \label{sec:amr_mark}

\subsubsection{Mesh optimality criteria: sub-module OPTIMALITY} \label{sec:amr_mark_optimality}

The strategy for selecting the elements to be refined is basically
based on a qualitative use of the estimated element-wise error
field. Qualitative approaches involve an empirical detection of the
elements to be refined. The best-known are quantile marking, Dörfler
marking or maximum marking, see~\cite{Dorfler-1996} for further
details. In this case, the optimal choice of the required parameters is not obvious,
since these parameters are not directly linked to a desired accuracy and may depend
on the problem under consideration, see~\cite{Barbie-2015,Ramiere-2019} for example.

In our developments, we prefer to rely on quantitative marking approaches,
related to an user-prescribed accuracy. To
this end, a maximum permissible element-wise error, denoted  $\xi_{T_{i}}^{\max}$,  has to be
defined. Two mesh optimality criteria, focusing either on a global or
local error, are studied.

$\bullet$ \textbf{ZZ criterion: equidistribution of the specific error.}

In~\cite{Zienkiewicz-1987}, Zienkiewicz and Zhu proposed a mesh
optimality criterion based on principle of error
equidistribution:
\begin{equation} \xi_{T_{i}}^{\max , \text{ZZ}} = \text{const}
  \quad \forall T_{i} \in \Omega^{h} .
\label{eqn:equidistribution_error_principle}
\end{equation}
The global absolute estimated error on $\Omega^{h}$ is then defined by \begin{equation}
\xi_{\Omega}^{\max , \text{ZZ}} = e_{\Omega} \omega_{\Omega} = \sqrt{N_{E}}\, \xi_{T_{i}}^{\max , \text{ZZ}}
\label{eqn:global_absolute_estimated_error_rewrite}
\end{equation}
with $N_{E}$ the number of mesh elements, $e_{\Omega}$ the
user-prescribed accuracy and $\omega_{\Omega}$ the energy norm
associated to the global strain. This value is obtained by summing the
elementary contributions:
\begin{equation}
\omega_{\Omega} = \left( \sum_{T_{i} \in \Omega^{h}} \omega_{T_{i}}^{2} \right)^{1/2}
\label{eqn:norme_globale_zz}
\end{equation}
where the strain energy $\omega_{T_{i}}^2$ on element $T_{i}$,
introduced in~\cite{Zienkiewicz-1987}, is given by
\begin{equation}
\omega_{T_{i}}^2 = \int_{T_{i}} \boldsymbol{\sigma}^{h} : \boldsymbol{\varepsilon}^{h} d T_{i} + \int_{T_{i}} \left( \boldsymbol{\sigma}^{\star} - \boldsymbol{\sigma}^{h} \right) : \left( \boldsymbol{\varepsilon}^{\star} - \boldsymbol{\varepsilon}^{h} \right) d T_{i}  \; \forall \; T_{i} \in \Omega^{h}.
\label{eqn:norme_locale_zz_2}
\end{equation}

Thus, the maximum permissible absolute element-wise error is then defined by \begin{equation}
\xi_{T_{i}}^{\max , \text{ZZ}} = e_{\Omega} \dfrac{\omega_{\Omega}}{\sqrt{N_{E}}} .
\label{eqn:critere_erreur_amr_seuil_zz}
\end{equation}
Criterion~\eqref{eqn:critere_erreur_amr_seuil_zz} is the most
widespread mesh optimality criterion known to date. 

$\bullet$ \textbf{LOC criterion: local element-wise error.}

In~\cite{Ramiere-2019}, a local optimality criterion (LOC) is introduced
and has been shown in~\cite{Koliesnikova-2021} to be more reliable for local error control
 in the framework of structural mechanics
problems. In addition, this criterion has been evaluated
in~\cite{Liu-2017} for contact problems. The maximal permissible
element-wise error distribution is directly governed by the admissible
local element-wise error:
\begin{equation}
\xi_{T_{i}}^{\max , \text{LOC}} = e_{\Omega , \text{LOC}} \; \omega_{T_{i}} .
\label{eqn:critere_erreur_amr_seuil_loc}
\end{equation}
As proved in~\cite{Ramiere-2019}, this automatically implies the
relative global error to be bounded above by $e_{\Omega ,
  \text{LOC}}$.\\

In the MFEM software library, the \texttt{mfem::MeshOperator::ApplyImpl} function
dedicated to the elements marking for refinement has been enriched
with implementations required to evaluate the maximal permissible errors: calculations of $\omega_{T_{i}}$,
$\omega_{\Omega}$, and, a fortiori, $\xi_{T_{i}}^{\max , \text{ZZ}}$
and $\xi_{T_{i}}^{\max , \text{LOC}}$.

\subsubsection{Detect elements: sub-module DETECT} \label{sec:amr_mark_detect}

From the maximum permissible element-wise error, the set of marked
elements is simply defined by
\begin{equation}
\mathcal{M}^{0} = \left\lbrace T_{i} \in \Omega^{h} \; / \; \xi_{T_{i}} > \xi_{T_{i}}^{\max} \right\rbrace .
\label{eqn:critere_detection_ensemble_elements_amr}
\end{equation}

In order to garantee the refined contact boundaries to remain
matched (required for the node-to-node contact pairing), the
$\mathcal{M}^{0}$ set has to be
enlarged. Elements located on the potential contact
boundaries whose paired element is marked for refinement, are then also
marked for refinement as proposed
in~\cite{Liu-2017}.
\begin{equation} \mathcal{M} =
  \mathcal{M}^{0} \cup \left\lbrace T_{j} \in \Omega^{h} \; / \; T_{i}
    \in \mathcal{M}^{0} \; \text{paired with} \; T_{j}
  \right\rbrace .
\label{eqn:critere_elargi_detection_ensemble_elements_amr}
\end{equation}
Here, we define finite elements as paired when all their nodes located
on contact boundaries are each other paired.

To that end, the MFEM \texttt{mfem::MeshOperator::ApplyImpl} function has been enriched
to receive paired elements as input.

\subsection{Stopping criteria: module STOP} \label{sec:amr_stop}

The AMR process obviously stops when the
set of marked elements becomes empty:
\begin{equation}
\mathcal{M} = \emptyset .
\label{eqn:critere_stop_ideal}
\end{equation}

However, for problems with singularities, such as contact mechanics
ones, criterion~\eqref{eqn:critere_stop_ideal} may never be reached
and infinite refinement in very localized zones may be performed. Other stopping criteria are therefore commonly added, allowing the refinement process to be ended once the mesh is considered as acceptable.

In this study, we may be interested in verifying either global or local accuracies. The refinement process is therefore analyzed using two stopping criteria. \\

$\bullet$ \textbf{Global stopping criterion.}

The stopping criterion most commonly used in the literature is a global
criterion based on the estimated absolute global error, see Eq.~\eqref{eqn:stop_acceptable_critere}.

\begin{equation}
\dfrac{\xi_{\Omega}}{e_{\Omega} \omega_{\Omega}} \leq 1 .
	\label{eqn:stop_acceptable_critere}
\end{equation}

$\bullet$ \textbf{Local stopping criterion.}

When the control of the local error is of interest (see detection
criterion~\eqref{eqn:critere_erreur_amr_seuil_loc}), a local stopping criterion, introduced
in~\cite{Ramiere-2019} and based on geometric considerations, is
used to stop the refinement process. Indeed, the respect of the
global stopping criterion~\eqref{eqn:stop_acceptable_critere} does not
guarantee any precision at local level. The local stopping criterion,
see Eq.~\eqref{eqn:critere_stop_acceptable_local},
ensures that the measure of the set of detected elements
$\mathcal{M}$ is below
a proportion ($\delta$) of the total measure of the domain.
\begin{equation}
\mu \left( \Omega_{\mathcal{M}} \right) \leq \delta \cdot \mu
\left( \Omega^{h} \right) 
\label{eqn:critere_stop_acceptable_local}
\end{equation}
with $\mu \left( \Omega_{\mathcal{M}} \right)$ the measure of
detected area $\overline{\Omega_{\mathcal{M}}} = \left\lbrace \cup \overline{T_{i}}
  \; / \; T_{i} \in \mathcal{M} \right\rbrace$, $\mu
\left( \Omega^{h} \right)$ the measure of global mesh, and $\delta$ the
user-defined local error control parameter.\\
With this stopping criterion, the obtained solution satisfies the
local optimality criterion~\eqref{eqn:critere_erreur_amr_seuil_loc}
over $\left( 1 - \delta \right) \times 100$ percent of the
domain. This local stopping criterion therefore has the advantage of automatically determining a discrete approximation of the critical area in the presence of contact singularities, which is of great interest.

\textbf{Remark:} Setting $\delta = 0$ makes the criterion equal
to~\eqref{eqn:critere_stop_ideal}. In addition, choosing $\delta \ll
1$ prevents a rough representation of the critical area. 

\subsection{Mesh refinement: module REFINE} \label{sec:amr_refine}

The mesh refinement is based on an isotropic hierarchical refinement,
consisting in dividing each marked quad/hexa element (also called
parent element) in $2^m$ (with $m$ the space dimension) sub-elements
of the same type (also called child elements), based on splitting each
side into $2$ equal parts. The refinement ratio is chosen to be equal to $2$
according to the study conducted in~\cite{Koliesnikova-2021}.

In addition, the 'single
irregularity rule' (also called 'single hanging node rule' or '2:1
balance rule'), argued in \cite{Nochetto-2012,Koliesnikova-2021} for example and widely used in
the literature, is applied. This has the effect of limiting to one the number
of irregular nodes on each edge. 

\begin{figure}[!h]
  \centering
	\includegraphics[scale=0.2]{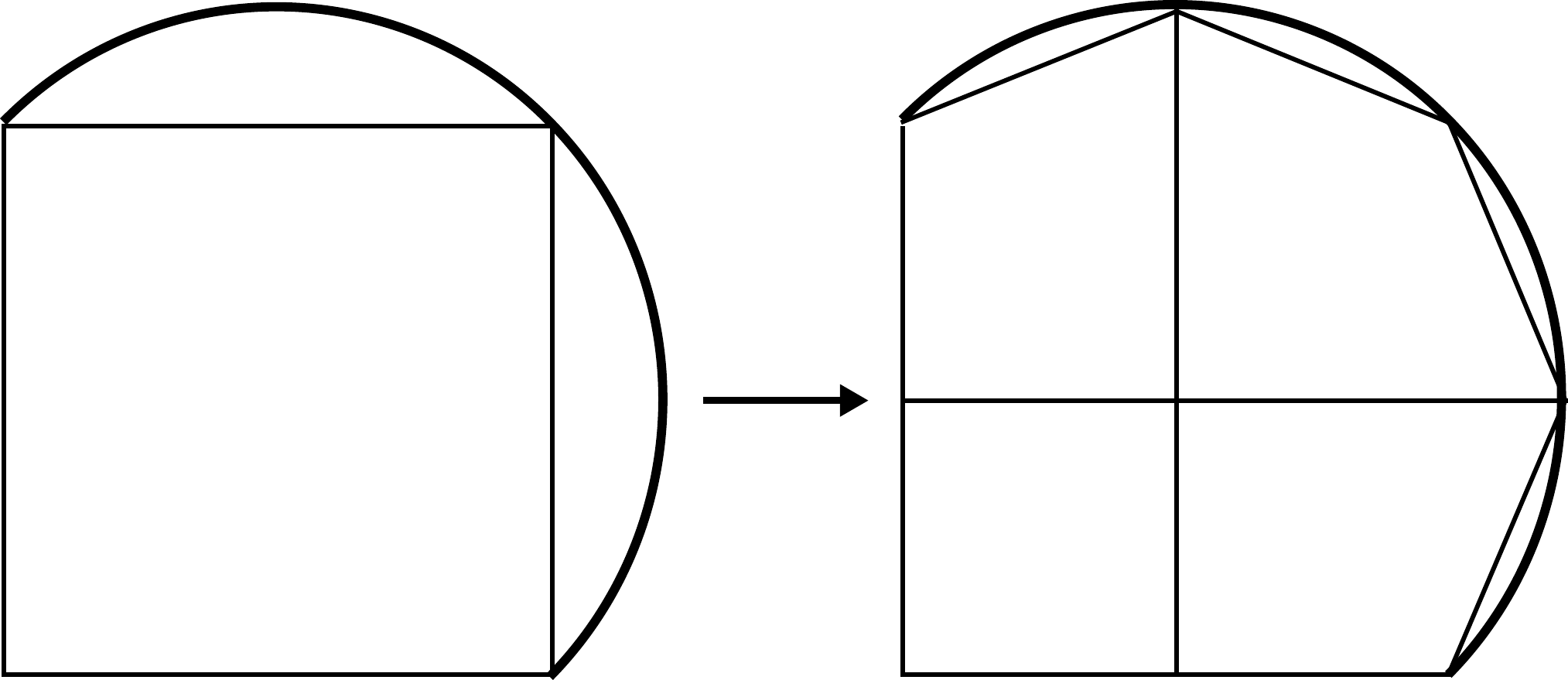}
	\caption{Quasi-hierarchical isotropic mesh refinement.}
	\label{fig:quasi-hierarchical-raff}
\end{figure}

However for first-order finite elements, this hierarchical
mesh refinement limits the geometry representation when considering
curved boundary, which is generally the case of contact boundaries.
In the literature (see \cite{Liu-2017,Koliesnikova-2021} for example), hand-made
'quasi-hierarchical' mesh refinement procedures have been
proposed. They mainly consist of updating the discretization of the
curved boundary by moving the nodes added by refinement onto the exact
curved geometry, see Figure~\ref{fig:quasi-hierarchical-raff}.
This procedure gives very satisfactory results, but
is intrusive and not generic as it requires knowledge of the
closed-form equation of the curved boundary during the refinement step.

\begin{figure}[!h]\centering
\subfloat[Initial coarse
mesh]{ \includegraphics[scale=0.2]{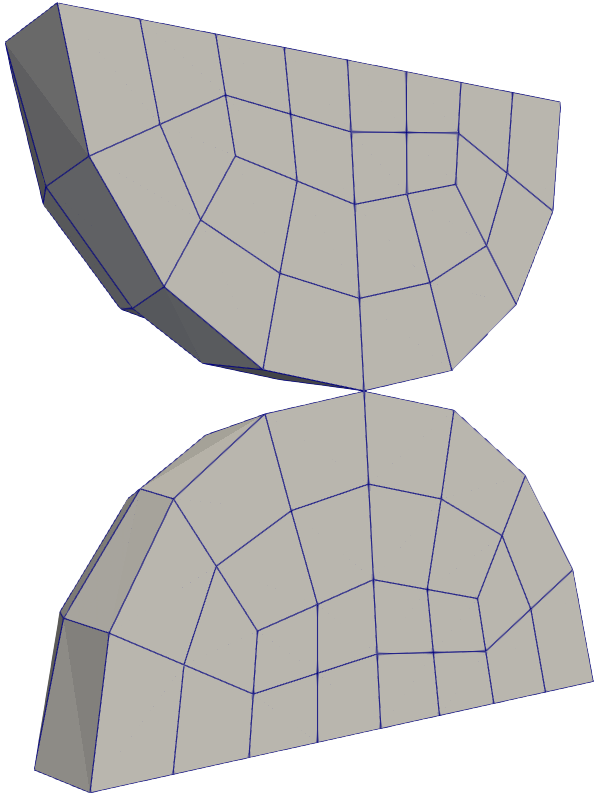} \label{fig:3d_initial_mesh}}
\hspace*{0.15cm}
\subfloat[Hierarchical AMR with isoparametric elements (first-order
transformation)]{ \includegraphics[scale=0.2]{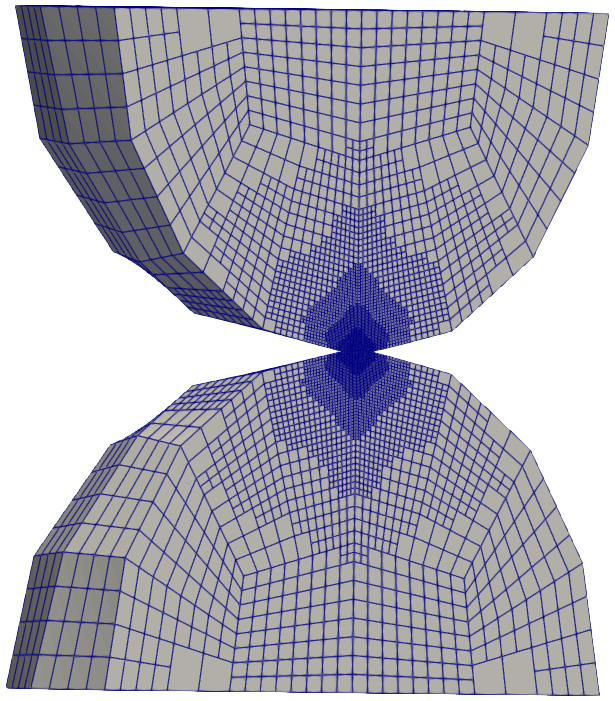} }
\hspace*{0.15cm}
\subfloat[Hierarchical AMR with super-parametric elements (6th-order transformation)]{ \includegraphics[scale=0.2]{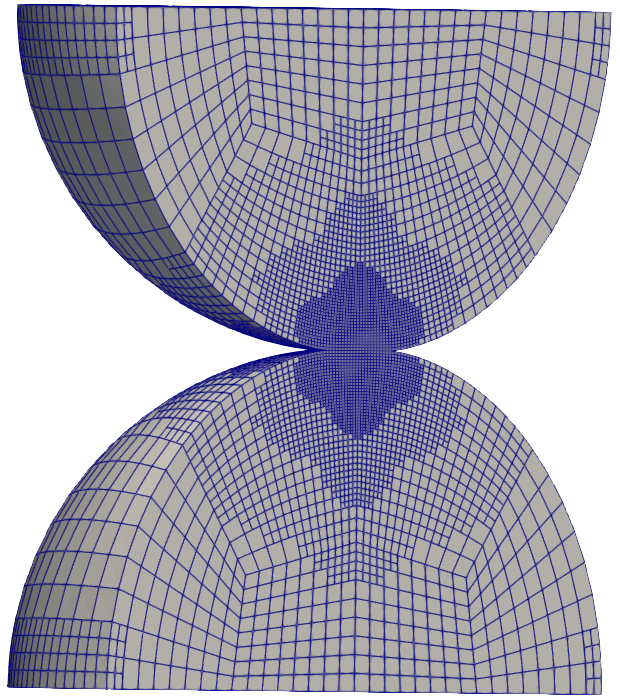} }
\caption{Hierarchical AMR for curved-edge contact problem with
  first-order finite element solution: isoparametric versus
  super-parametric elements.}
\label{fig:hierarchical-raff-superparametric}
\end{figure}

In this contribution, shape preservation of hierarchically refined
geometries is automatically guaranteed by the use of super-parametric
finite elements.  These elements
are based on a geometric transformation of higher degree of
interpolation than the solution's basis functions.
With this type of element, the surface approximation error during
refinement can be very low compared to the one obtained with a first-order
transformation, see Figure~\ref{fig:hierarchical-raff-superparametric}. These elements are therefore ideally suited to contact
mechanics problems where the active contact zone must be charactized
precisely.

This type of super-parametric element is defined during the mesh design
phase and is properly supported by the MFEM software environment. As a
result, the discrete geometry does not need to be updated
during the refinement process.

The MFEM refinement procedure, named \texttt{mfem::MeshOperator::Apply}, is invoked to
carry out the chosen non-conforming h-adaptive mesh refinement
strategy~\cite{Cerveny-2019, Anderson-2020}. The single irregularity
rule is activated by setting the \texttt{nc\_limit} parameter to $1$.

\subsection{Non-conforming solution: module SOLVE} \label{sec:non_conforming_mfem}

\begin{figure}[!h]
  \centering
\includegraphics[scale=0.2]{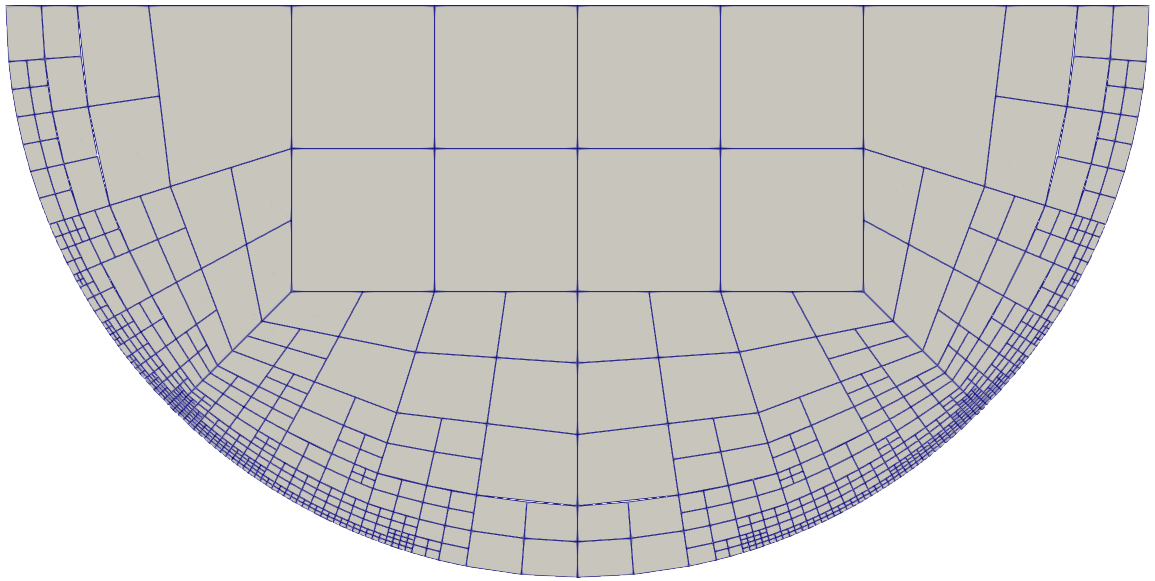}
\caption{Example of a non-conforming h-adaptive refined mesh.}
\label{fig:maillage_non_conforme}
\end{figure}

Refinement of a quadrangular or hexahedral elements mesh using the
hierarchical h-adaptive method results in the appearance of hanging
nodes, see Figure~\ref{fig:maillage_non_conforme}. Thus, the
non-conforming degrees of freedom (DOFs) must be constrained in order
to obtain a conforming solution.

The major task in the use of this type of refinement method is to
build the linear system to be inverted. The choice made in MFEM and
detailed in~\cite{Cerveny-2019} is to define (and solve) a problem
restricted to conforming DOFs only. Its principle is recalled here after in
order to present its extension to contact problems.

The global stiffness matrix (bilinear form) $\tilde{\textbf{K}}$ and
the load vector of the linear finite element system
$\tilde{\textbf{L}}^{ext}$ are assembled taking into account all the
mesh nodes, without worrying about non-conformities. Thus, the
resulting linear system writes
\begin{equation} \tilde{\textbf{K}}
  \tilde{\textbf{U}} = \tilde{\textbf{L}}^{ext}
\label{eqn:assemblage_global}
\end{equation}
with $\tilde{\textbf{U}}$ the vector containing all DOFs.\\
In order to restrict this system to conforming DOFs, a prolongation
operator $\textbf{P}$ is built:
 \begin{equation}
\textbf{P} = \begin{pmatrix}
\textbf{Id} \\ 
\textbf{W}
\end{pmatrix}
\label{eqn:prolongation_operator}
\end{equation}
with $\textbf{Id}$ the identity matrix and $\textbf{W}$ the so-called
interpolation matrix. This interpolation matrix $\textbf{W}$ expresses
the non-conforming DOFs to be constrained as linear combinations of
the conforming DOFs using the finite element basis functions
(see~\cite{Cerveny-2019, Anderson-2020} for further details on the
procedure for building these matrices).
Hence, we have
\begin{equation}
\tilde{\textbf{U}} = \textbf{P} \textbf{U} 
\label{eqn:prolongation_form}
\end{equation}
with $\textbf{U}$ the conforming DOFs. Similarly, the test functions are
projected onto the conforming space thanks to the prolongation
operator $\textbf{P}$. Then, the linear
system~\eqref{eqn:assemblage_global} becomes
\begin{equation}
  \label{eq:systeme-projete}
\textbf{P}^{T} \tilde{\textbf{K}} \textbf{P} \textbf{U} = \textbf{P}^{T} \tilde{\textbf{L}}^{ext} .
\end{equation}
By defining the conforming stiffness matrix
\begin{equation}
\textbf{K} = \textbf{P}^{T} \tilde{\textbf{K}} \textbf{P} ,
\label{eqn:bilinear_form}
\end{equation}
and the conforming load vector
\begin{equation}
\textbf{L}^{ext} = \textbf{P}^{T} \tilde{\textbf{L}}^{ext},
\label{eqn:linear_form}
\end{equation}
the linear system to be solved to obtain the solution $\textbf{U}$ on conforming DOFs writes
 \begin{equation}
\textbf{K} \textbf{U} = \textbf{L}^{ext}.
\label{eqn:assemblage_conforme}
\end{equation}

Once the finite element solution $\textbf{U}$ is found, it is extended
to all DOFs of the mesh through the prolongation
operation~\eqref{eqn:prolongation_form}.\\


The same projection
operations need to be performed on the contact terms appearing in the
penalized linear
system~\eqref{eqn:mechanical_equilibrium_matrix_system}.
As said in Section~\ref{sec:discrete-form}, we have chosen to rely on
the \texttt{mfem::PenaltyConstrainedSolver} class of MFEM which works on
matrices and vectors associated to conforming DOFs. To this end, it is
sufficient to define the active contact pairing matrix on conforming
DOFs by
 with 
\begin{equation}
\hat{\textbf{B}} = \tilde{\textbf{B}}[\mathcal{A}] \textbf{P}
\label{eqn:contact_terms}
\end{equation}
where $\tilde{\textbf{B}}$ is the potential contact pairing matrix
assembled taking account all the DOFs, and $\mathcal{A}$ the set of
active contact nodes, see Section~\ref{sec:discrete-form}.\\
Then, the contact stiffness matrix built by the
\texttt{mfem::PenaltyConstrainedSolver} class writes
 \begin{equation}
\hat{\textbf{K}}_{C} = \textbf{P}^{T}  \tilde{\textbf{B}}[\mathcal{A}]^{T}\tilde{\textbf{B}}[\mathcal{A}] \textbf{P},
\label{eqn:contact_matrix}
\end{equation}
which correspond to the projection to conforming DOFs of
the global contact stiffness matrix $\tilde{\textbf{K}}_{C}=
\tilde{\textbf{B}}[\mathcal{A}]^{T}\tilde{\textbf{B}}[\mathcal{A}]$.
The projection of the contact right-hand term is also consistent:
\begin{equation}
\hat{\textbf{B}}^{T} \hat{\textbf{D}} = \textbf{P}^{T} \tilde{\textbf{B}}[\mathcal{A}]^{T} \hat{\textbf{D}} .
\label{eqn:contact_linear_form}
\end{equation}
Hence the contact right hand term $\hat{\textbf{D}}$ given of input of
\texttt{PenaltyConstrainedSolver} is directly
\begin{equation}
  \hat{\textbf{D}} = \tilde{\textbf{D}}[\mathcal{A}] .
\end{equation}

It should be noted that these algebraic operations are not directly available in
MFEM classes and have been implemented specifically.

\subsection{Algorithmic treatment of the AMR-Contact solution process} \label{sec:amr_contact_algorithm}

The proposed AMR-Contact solution process is ruled by two nested
iterative loops, described in
Algorithm~\ref{alg:amr_process_algo}. This algorithm calls
Algorithm~\ref{alg:contact_solution_process_algo} in which the
definition of the initial set of active contact nodes is generally
performed considering the solution at the
previous AMR iteration (for computational efficiency).

\begin{algorithm}[!h]
\caption{Combined AMR-Contact solution process}
\label{alg:amr_process_algo}
\begin{algorithmic} 
\STATE \textbf{Input: Data of the problem, maximum number of AMR loop iterations $n_{max}$}
\STATE \textbf{Output: Finite element solution on non-conforming refined mesh}
\STATE \textbf{Initialization}
\STATE \hspace*{0.25cm} $n \leftarrow 1$, $\textit{\textbf{stop}}_{\textit{\textbf{AMR}}} \leftarrow false$
\WHILE{$\textit{\textbf{stop}}_{\textit{\textbf{AMR}}} = false$ and $n \leq n_{max}$}
\STATE Assemble the bilinear and linear forms
\STATE \textbf{Restriction step.}
\STATE \hspace*{0.25cm} Form the conforming linear system
\STATE Perform the contact solution process (SOLVE): call Algorithm~\ref{alg:contact_solution_process_algo}
\STATE \textbf{Prolongation step.}
\STATE \hspace*{0.25cm} Reconstruction of the finite element solution
on all DOFs
\STATE Detect elements to be refined (ESTIMATE-MARK)
\STATE Update $\textit{\textbf{stop}}_{\textit{\textbf{AMR}}}$ with
the stopping criterion (STOP)
\IF{$\textit{\textbf{stop}}_{\textit{\textbf{AMR}}} = false$}
\STATE Refine detected elements (REFINE):
\STATE \hspace*{0.25cm} Create the new non-conforming refined mesh
\STATE Rebalance the mesh partitioning (see Section
\ref{sec:mesh_partitioning})
\STATE $n \leftarrow n + 1$
\ENDIF
\ENDWHILE
\end{algorithmic}
\end{algorithm}

\section{Parallel strategy} \label{sec:parallel_strategy}

Since accurate prediction of contact mechanics problems with mesh
adaptation generally requires a large number of elements and a large
sparse system to solve, distributing tasks among processes is an important matter. With the need for parallel codes becoming increasingly important, algorithms have been
implemented in finite element libraries such as MFEM for handling
large non-conforming meshes~\cite{Anderson-2020}. However,
adapting them for addressing contact mechanics problems remains a
challenging task. In this section, we describe a parallel solution
that is dedicated to this type of problem. To this end, the most common distributed-memory
parallel processing paradigm is considered, namely Message Passing
Interface (MPI).

\subsection{MFEM non-conforming mesh partitioning} \label{sec:mfem_non_conforming_partitioning}

An effective combination of contact mechanics and non-conforming
h-refinement strategies in a distributed-memory parallel processing
framework is based, first and foremost, on a relevant mesh partitioning. In the
MFEM library, the standard type of partitioning is element-based,
i.e. each element is assigned to one of the MPI process
(see~\cite{Cerveny-2019} for further details).  The set of (finite)
mesh elements is divided into $R$ disjoint regions, where $R$ is the
number of MPI tasks. The data related to vertices, edges and faces at
the boundary of each task's region are duplicated by adding ghost
elements (copy of foreign elements that are owned by other MPI processes), so that each region can be treated as an isolated submesh
and processed in parallel. The prolongation operator $\textbf{P}$ (see
Eq.~\eqref{eqn:prolongation_operator}) is extended to
the processing of ghost DOFs so that the solution of the linear system, obtained after the solver step in the SOLVE module, concerns the non-ghost conforming DOFs only.

By default in the MFEM library, when the mesh is hierarchically
refined, child elements are preassigned to their parent MPI task. In
case of non-conforming h-adaptive refinement, imbalance between MPI
processes occurs, and a rebalance of elements is worthwhile. The \texttt{mfem::ParMesh::Rebalance} function can be called
to produce a new mesh partition. Used as such (without any argument),
it corresponds to a black-box function that builds a
partitioning solution based on the classic geometric Space-Filling Curve
(SFC) method~\cite{Pilkington-1994, Aluru-1997} via Hilbert curves.

In the MFEM library, the embedded functions for creating mesh partitioning are
not currently suited to contact mechanics problems. Indeed, if they
are applied as such, the paired contact nodes have no guarantee to be on the
same processes. Thus, the construction of the contact pairing matrix $\textbf{B}$ would
require, at least, extra MPI communications. In this paper, we
present a specific strategy that provides a new partitioning function. Adapting the mesh partitioning functions proposed by MFEM
to contact problems is a future task. On the other hand, MFEM's other
main features (assembly, linear system formation, solving, etc.) are
still used.

%
%

\subsection{Contact-based mesh partitioning algorithm} \label{sec:mesh_partitioning}

The dedicated algorithm proposed here guarantees a mesh partitioning
whose the contact paired nodes are on the same processes during the
refinement iterations.
Practically, this involves
that the construction of the local contact pairing matrix $\textbf{B}$ does
not depend on communications with other processes. Hence, contact
stiffness matrices and load contact vectors are assembled locally on
each MPI task, as if each region were treated as an isolated mesh by
the sequential part of the code. This highly
reduces the amount of inter-process communications during the
refinement iterations, enhancing code reachability and parallel
scalability. \\
\textbf{Remark:} to simplify the code writing, elements lying on 
potential contact boundaries are not handled by the same MPI processes than
non-contact elements of the solids.

The partitioning algorithm is implemented in the sequential part of the
code, i.e. without any communication between the MPI tasks.
This algorithm is responsible for assigning each element to the MPI task to which it will belong once the rebalancing is carried out. Here, the
proposed algorithm provides a mesh partition to the \texttt{mfem::ParNCMesh::Rebalance}
function, which then accordingly redistributes the elements between the MPI
tasks. The algorithm is designed to ensure that, once the mesh has
been redistributed, the work load is as equal as possible between the processes handling the contact elements on the one hand, and
the other tasks on the other hand.

After applying the AMR algorithm (and before
rebalancing), MPI regions initially handling contact elements can
contain contact and non-contact elements.
Generally speaking, in our partitioning algorithm, the contact and non-contact elements of
the MPI regions are treated separately. The formula
\begin{equation}
R_{C} = \left \lceil c \dfrac{N_{E, C}}{N_{E}} R \right \rceil
\label{eqn:nb_contact_procs}
\end{equation}
defines the number of MPI tasks processing contact elements. In this
equation, $N_{E, C}$ represents the (even) number of contact elements,
$N_{E}$ the total number of mesh elements (both after AMR iteration), $R$ the
total number of MPI processes, $c > 0$ an a priori given
coefficient and $\left \lceil x \right \rceil$ the ceiling function of
the scalar $x$ (smallest integer greater or equal to $x$). In this
formula, choosing $c > 1$ will lead to have fewer elements on average
on MPI tasks processing contact elements than on the rest of the MPI
processes, and inversely for $c < 1$. In addition, $R_{C}$ is
reassessed at each mesh refinement iteration as $N_{E, C}$ and
$N_{E}$ may have changed.

The result obtained with Equation~\eqref{eqn:nb_contact_procs} must be corrected in the following cases: \begin{equation}
\begin{cases}
R_{C} = \dfrac{N_{E, C}}{2} & \text{ if } R_{C} > \dfrac{N_{E, C}}{2} \\ 
R_{C} = R - 1  & \text{ if } R_{C} \geq R \; \text{and} \; N_{E} > N_{E, C} .
\end{cases}
\label{eqn:nb_contact_procs_correction}
\end{equation}
The first criterion guarantees that there is at least one pair of
contact elements per MPI task processing contact elements. The second
criterion ensures that at least one MPI process is retained for
processing non-contact elements. Obviously, if there are contact
elements only, the second criterion does not occur.

Then, for each region $r$ containing $N_{E,r}$ elements, the
$N_{E,C,r}$ contact elements (with $N_{E,C,r}=0$ in non-contact MPI regions) must be assigned to
the $R_{C}$ first MPI processes handling the contact elements, while
the $(N_{E,r}$ - $N_{E,C,r})$ other elements must be assigned to
the following $(R - R_{C})$ MPI tasks. 


As the main feature of our mesh partitioning algorithm is to ensure
paired nodes to be on the same MPI task, a special attention is
devoted to the distribution of the contact elements.  Thanks to the
node-to-node pairing algorithm, the element-based mesh partitioning
algorithm can work on contact paired elements (see definition in
Section~\ref{sec:amr_mark_detect}).  Practically, in each MPI region
$r$, the contact elements of the first solid are distributed (quite)
equally according to their local number on the $R_{C}$ first MPI
processes. These elements are assigned locally to MPI processes via a
strategy close to SFC method applied to contact elements
only. However, it differs from this method in two notable
points. Firstly, unlike the SFC method via Hilbert curve, which performs a global
renumbering of elements in order to partition them, thereby inducing
communications between MPI tasks, only the local numbering of elements
is taken into account here for partitioning. Thus, no communication
between MPI processes is necessary. Secondly, as for a region $r$, there
could be fewer contact elements than MPI tasks processing the contact
elements ($N_{E,C,r} / 2 < R_{C}$), the number of the starting
process for partitioning is defined as $r\;\mathrm{mod}\;{R_C}$
(with 'mod' the modulo) to avoid potential imbalances. The number of
the ending process is then
$(r+N_{E,C,r} / 2)\;\mathrm{mod}\;{R_C}$.
Then, the contact elements in the second solid
are affected to the same MPI process than their paired element in the
first solid.

Then, the remaining (non-contact) elements must be distributed in a
balance way on the remaining subset of MPI tasks. This is done
according to a routine based on element equidistribution similar to
that applied for the contact elements of the first solid.

Finally, once this new partitioning has been produced, the MFEM \texttt{mfem::ParNCMesh::RedistributeElements} function, called within the \texttt{mfem::ParNCMesh::Rebalance} function, handles transferring elements
between MPI tasks, while respecting the specified partitioning. In addition, in the
proposed AMR-Contact solution strategy, the rebalancing step is
performed at the end of each AMR iteration, as this is a low-cost
step~\cite{Cerveny-2019}. \\

\begin{figure}[!h]\centering
\includegraphics[scale=0.4]{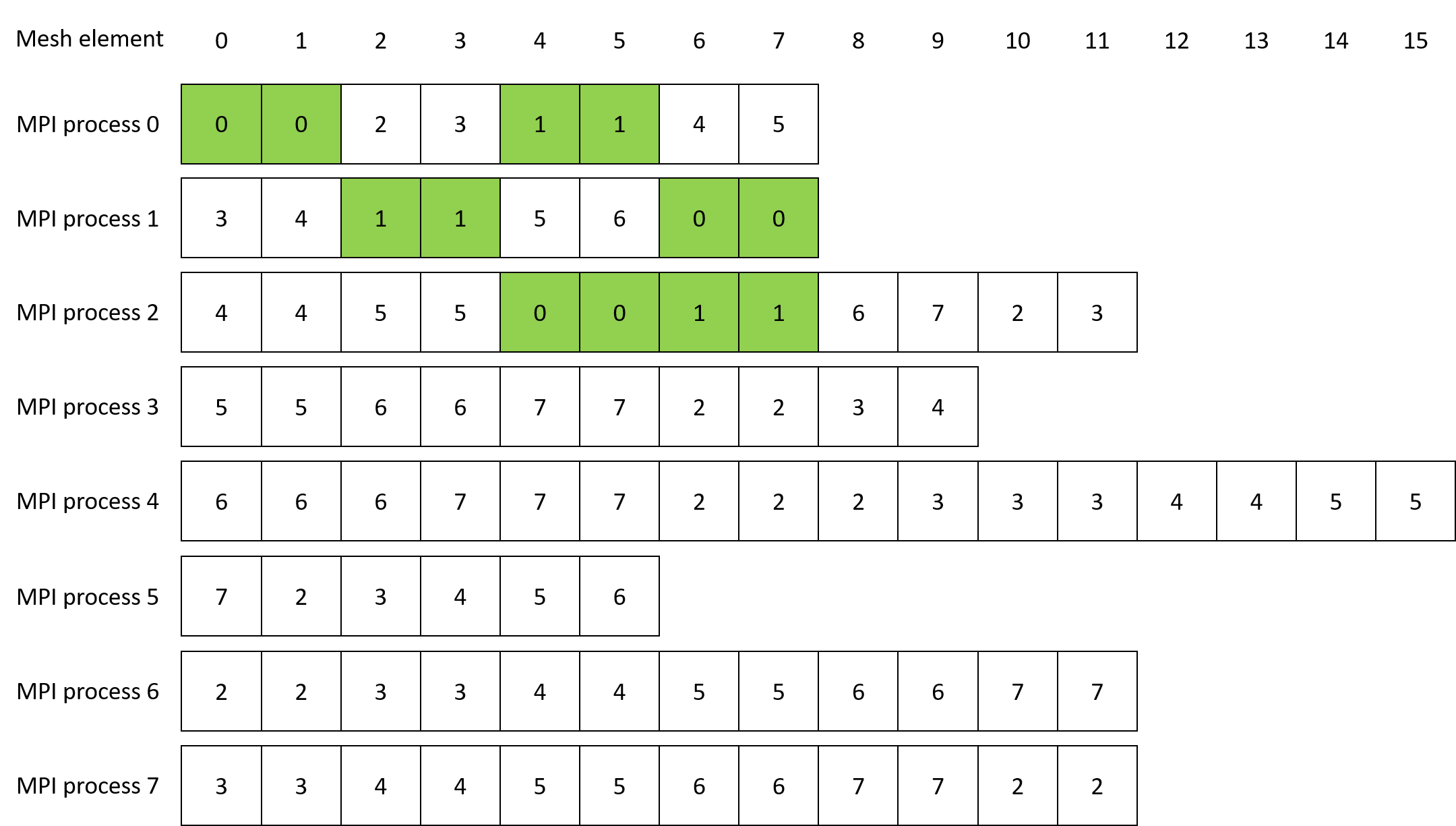}
\caption{Illustration of the partitioning principle for 8 MPI processes. In green = 'super-elements', in white = non-contact elements. The number of boxes per process corresponds to the number of 'super-elements' and non-contact elements that the processes have after applying AMR and before rebalancing the elements. The number in each cell specifies the process number to which the contact 'super-element' or the non-contact element will be assigned.}
\label{fig:equidistributed_partitioning_illustration}
\end{figure}

An example of a partitioning produced by our strategy is shown in
Figure~\ref{fig:equidistributed_partitioning_illustration}. In this example, a local mesh refinement step has been performed, and the elements (represented here by boxes) must be partitioned across 8 MPI tasks, which initially have a variable number of elements. At the current AMR iteration, there are 3 MPI tasks handling contact elements, which, after application of the refinement step, contain both contact elements and non-contact elements. In Figure~\ref{fig:equidistributed_partitioning_illustration}, contact elements are grouped by paired contact elements, and only the local number of the contact element of the first solid is considered. In this way, a contact element of the second solid is no longer identified by its local number, but by the local number of the contact element of the first solid to which it is paired. Each grouping of contact elements characterizes a 'super-element' and is illustrated in green in Figure~\ref{fig:equidistributed_partitioning_illustration}. Non-contact elements are shown in white. \\

In addition, at the current AMR iteration, there are also 5 MPI tasks
handling non-contact elements, which, once the refinement step has
been completed, contain obviously non-contact elements only.  For each
MPI process, the elements (or pairs of contact elements) are numbered
locally, and their local number is defined at the top row of
Figure~\ref{fig:equidistributed_partitioning_illustration}.  After
applying Equation~\ref{eqn:nb_contact_procs} and
criteria~\ref{eqn:nb_contact_procs_correction}, $R_{C}$ will be (in
this example) 2 at the next AMR iteration, and $R - R_{C}$ will
therefore be 6.  The MPI tasks handling the contact elements will
therefore be processes 0 and 1, while the processes handling the
non-contact elements will be processes 2, 3, 4, 5, 6 and 7.  The shift in
the starting process number for partitioning is also highlighted with the partitioning example presented in
the figure.
Indeed, for each process, the process number assigned to the mesh element corresponding to the first green (respectively white) box starting from the left is given by $\text{r} \; \mathrm{mod} \; R_{C}$ (respectively $\text{r} + R_{C}$ if $\text{r} + R_{C} < R$ or $\text{r} + R_{C} \; \mathrm{mod} \; R - R_{C}$ otherwise, with 'mod' being the modulo operation): for example, the first white box of process 4 has been assigned the number 6 and that of process 5 has been assigned the number 7. Similarly, this figure allows us to observe that the end process number (for contact elements: $\left( \text{r} + N_{E,C,\text{r}} / 2 \right) \; \mathrm{mod} \; R_{C}$~; for non-contact elements: $\text{r} + R_{C} + N_{E,\bar{C},\text{r}}$ if $\text{r} + R_{C} + N_{E,\bar{C},\text{r}} < R$ or $\text{r} + R_{C} + N_{E,\bar{C},\text{r}} \; \mathrm{mod} \; R - R_{C}$ otherwise) differs depending on the process considered (and the number of elements it has to redistribute). For example, for the contact elements of process 1, the end process is 0, and for the non-contact elements of process 4, the end process is 5.

Furthermore, it should be noted that the case where $N_{E,C,r} / 2 < R_{C}$ (which would mean that the number of green boxes is strictly less than $R_{C}$) does not happen
in the example proposed here. Indeed, for each process with contact elements (initially
MPI rank = 0, 1, 2), the
number of green cells is greater than $R_{C} = 2$.

Then, the
proposed partitioning algorithm provides a new local partition of the
mesh which is that given in
Figure~\ref{fig:equidistributed_partitioning_illustration}, where the
number inscribed in each box indicates the MPI process to which the
element (or pair of contact elements) is newly assigned. \\

\begin{figure}[!h]\centering
\subfloat[][]{ \includegraphics[scale=0.25]{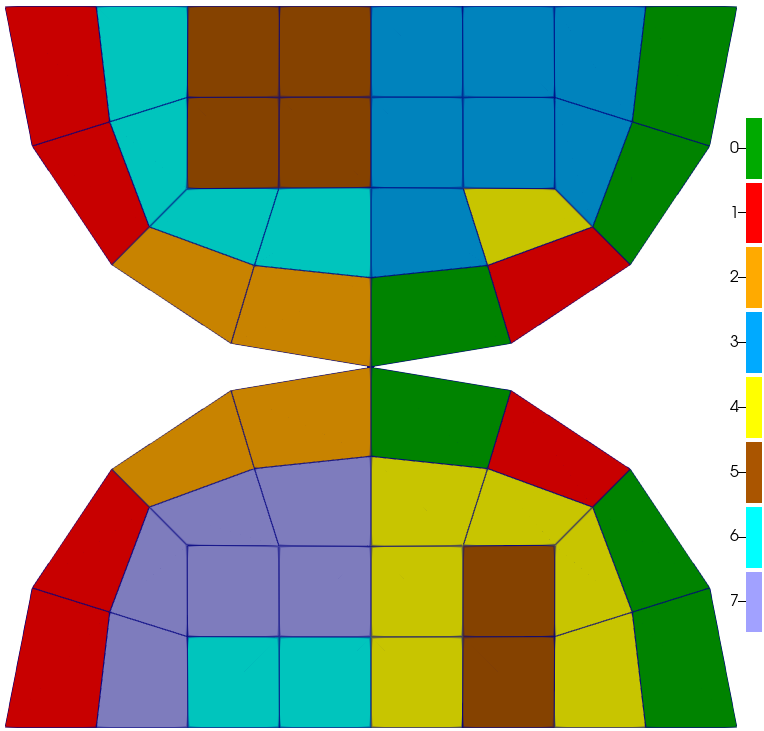} \label{fig:equidistributed_partitioning_amr0_2Dhertzdisques}}
\hspace*{1\baselineskip} 
\subfloat[][]{ \includegraphics[scale=0.25]{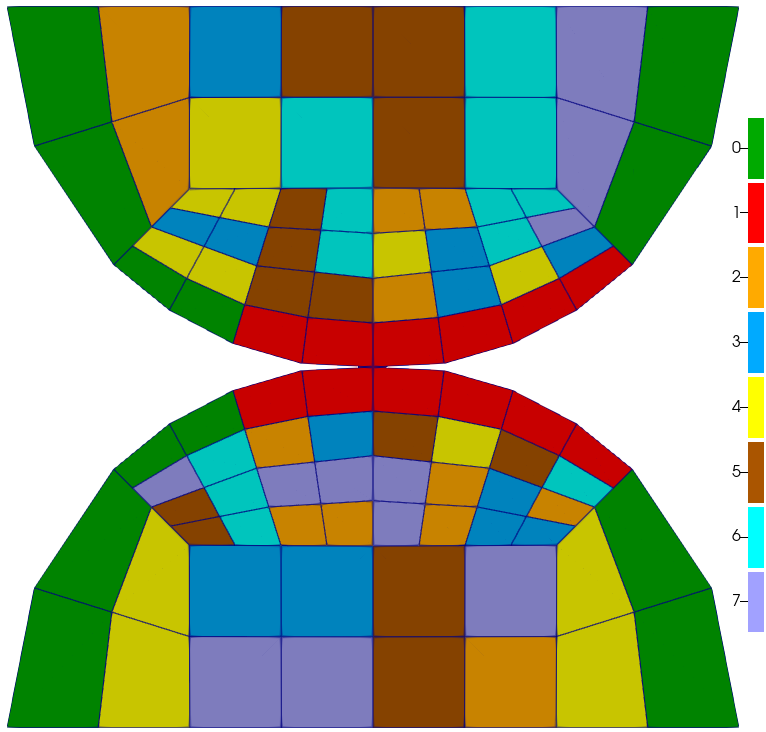} \label{fig:equidistributed_partitioning_amr1_2Dhertzdisques}}
\caption{Contact-based load balancing illustration for a 2D mesh,
  typical of the Hertzian problem - Each color identifies a specific MPI task - Mesh partitionings at first and second AMR iteration in (a) and (b) respectively.}
\label{fig:equidistributed_partitioning_amrfirsts_2Dhertzdisques}
\end{figure}

Applying the proposed strategy to the example presented in Figure~\ref{fig:equidistributed_partitioning_illustration} leads to the mesh partitioning illustrated on 2D meshes in Figure~\ref{fig:equidistributed_partitioning_amr1_2Dhertzdisques}. The initial partitioning is reported on Figure~\ref{fig:equidistributed_partitioning_amr0_2Dhertzdisques}\\

This figure displays the MPI task numbers assigned to the elements of
a 2D mesh used for an contact simulation with AMR (here, $c = 1$). The
calculation was performed using 8 cores to facilitate understanding of
the figure, and the first two AMR iterations are reported. One can notice on this figure
 that the initial global numbering of mesh elements (performed by
the Gmsh mesher~\cite{Geuzaine-2009}) has an influence on the mesh
partitioning, both in the first and subsequent AMR iterations. Indeed,
the global numbering of the mesh is used for initial partitioning and influences the new local numbering of the elements, obtained once the mesh partitioning is applied. Additionally, we should mention that, although the elements are staggered among the MPI processes, the MFEM library handles this well internally. MFEM manages the lack of element cohesion through ghost elements, which are automatically added and tracked around elements within a given MPI process. These ghost zones act as internal boundaries and are added to transparently handle MPI communications with neighboring MPI processes.\\


%
%

\begin{figure}[!h]\centering
\subfloat[][]{ \includegraphics[scale=0.2]{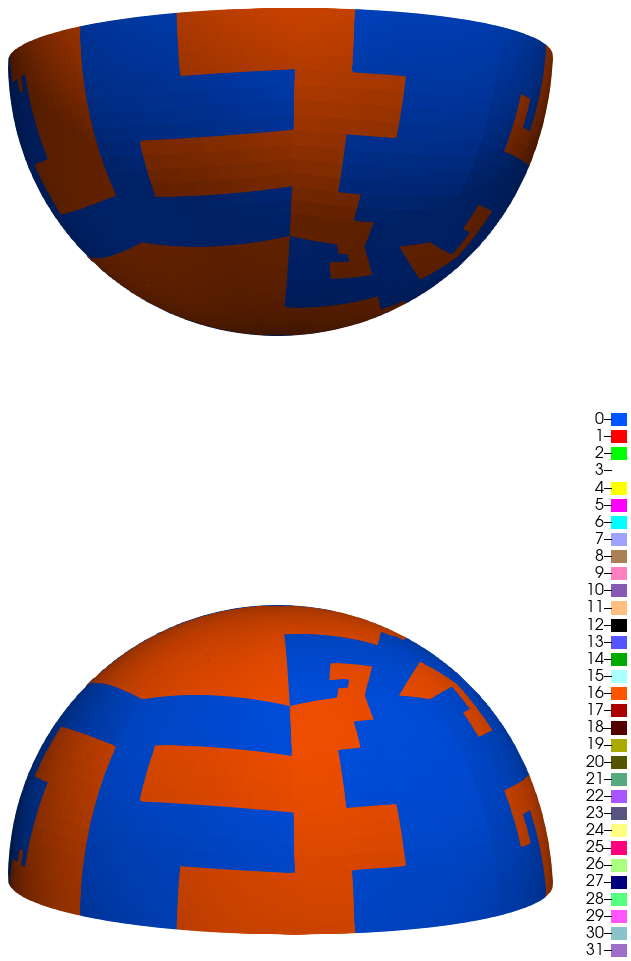} \label{fig:3d_amr_err2_glob_loadbalance2}}
\hspace*{1\baselineskip} 
\subfloat[][]{ \includegraphics[scale=0.175]{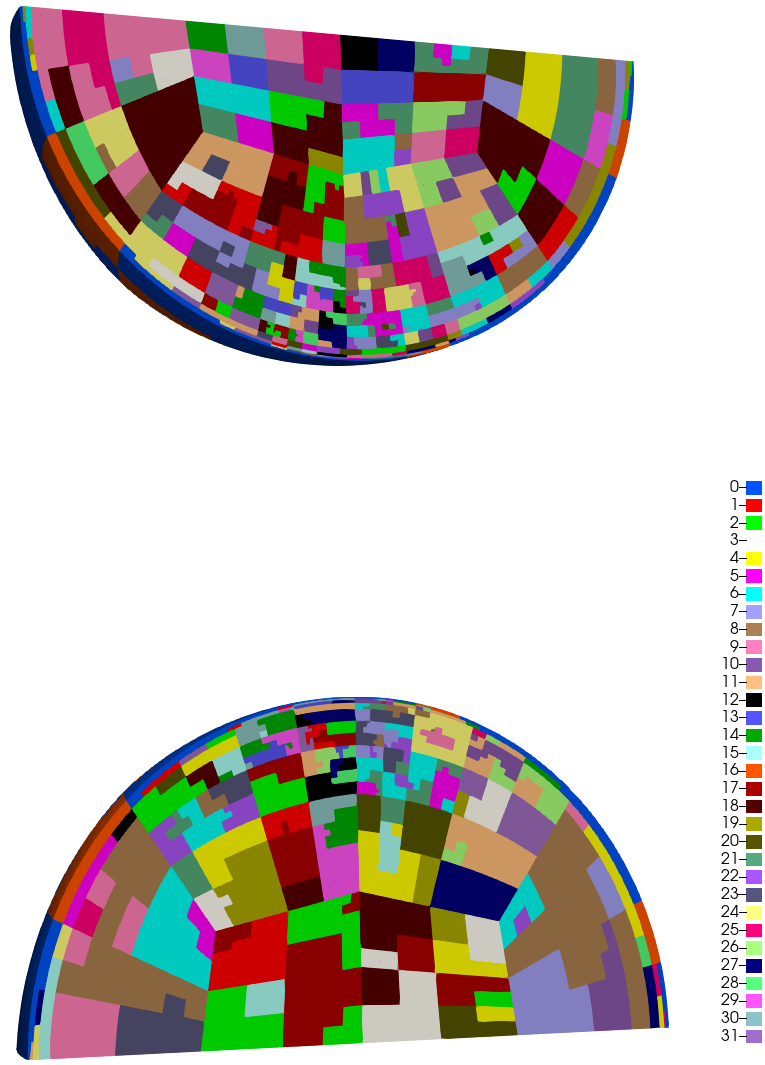} \label{fig:3d_amr_err2_glob_loadbalance3}}
\caption{Contact-based load balancing illustration for a 3D mesh,
  typical of the Hertzian problem - Color representing the MPI task.}
\label{fig:3d_amr_err2_glob_loadbalance}
\end{figure}

As mentioned earlier, this mesh partitioning method, which was set up to handle contact, ensures
the contact terms to be
built without inter-MPI process communications.  On the other
hand, the proposed mesh partitioning induces MPI regions composed of
elements that are scattered throughout the mesh. This can be seen even
better in
Figure~\ref{fig:3d_amr_err2_glob_loadbalance} which displays the MPI
task numbers assigned to the elements of a 3D mesh used for a
AMR-contact simulation (here, $c = 1$). The calculation was performed using 32
cores to have an acceptable number of colors, and the mesh edges are
not represented to better distinguish the color changes that means a
change of MPI task between two neighboring elements.

As expected contact paired nodes are on the same MPI task, but many
MPI-task interfaces are created inside each solids.

\section{Application to the Hertzian contact problem} \label{sec:hertzian_problem}

\begin{figure}[!h]\centering
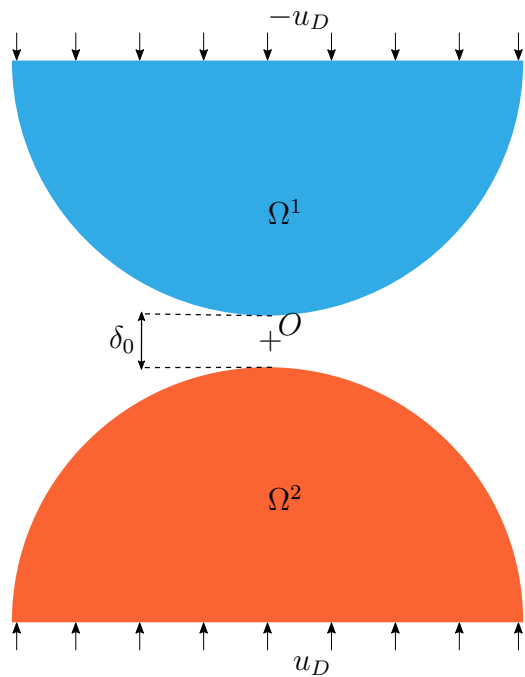
\caption{2D Hertzian contact problem.}
\label{fig:hertz_contact}
\end{figure}

This section is devoted to illustrating the efficiency of the proposed
parallel strategy and its suitability to contact mechanics problems
with curved edges. To this end, academic Hertzian benchmark problems
from~\cite{Hertz-1882} are considered is 2D and 3D. There are
well-known reference test cases in contact mechanics, especially as
analytical expressions of the stress field in the contact area are
available to verify numerical strategy accuracy. The mechanical
problem under consideration is defined
by~\eqref{elastostatic_contact_problem_hpp}.  The two 2D half-disks
(resp. two 3D half-spheres) brought into contact are
isotropic elastic materials with the same Young's modulus $E = 210$~GPa and a
Poisson's ratio $\nu = 0.3$. Their radius is equal to $R = 2$~m.
Dirichlet boundary
conditions of value $u_{D} = \dfrac{\delta_{0}}{2} + \alpha R$, where $\delta_{0} = 2$~m is the initial gap between solids and $\alpha$ a positive coefficient, are imposed on the planar
boundaries such as depicted for example in 2D
in Figure~\ref{fig:hertz_contact}.

Two
combinations of detection and stopping criteria for AMR are
considered:
\begin{itemize}
\item \textbf{AMR Combination 1:} ZZ optimality
criterion~\eqref{eqn:critere_erreur_amr_seuil_zz} and global stopping
criterion~\eqref{eqn:stop_acceptable_critere};
\item \textbf{AMR Combination 2:} LOC optimality
criterion~\eqref{eqn:critere_erreur_amr_seuil_loc} and local stopping
criterion~\eqref{eqn:critere_stop_acceptable_local}.
\end{itemize}

The quantities of interest in our simulations are the relative global estimated error in energy norm
$\gamma_{\Omega} = \xi_{\Omega} / \omega_{\Omega}$,
the measure of the zone still marked for refinement at the end of the
AMR process evaluated by the ratio
$\eta = \mu \left( \Omega_{\mathcal{M}} \right) / \mu
\left( \Omega^{h} \right)$, the number of mesh elements
$N_{E}$ and the number of DOFs $N$. These quantities of interest make it possible to quantitatively compare the two considered AMR combinations.

\subsection{Justification of the penalty parameter $k_{N}$}

Although the penalty method is a classical and widespread method for
numerical treatment of contact problems, see for
example~\cite{Kikuchi-1988, Laursen-2002, Wriggers-2006}, the choice
of an optimal penalty parameter (see $k_N$ is Eq.~\eqref{eqn:force-contact-penalization}) remains the main limitation of this
method. Indeed, the generated interpenetration error is directly
related to the penalization parameter $k_{N}$, and goes to zero only
in an asymptotic case, when $k_{N}$ tends towards infinity. However, a
too large penalty coefficient consistently produces ill-conditioning
of the stiffness matrix $\hat{\textbf{K}}_{C}$, leading to severe
convergence issues in the solution of
Problem~\eqref{eqn:mechanical_equilibrium_matrix_system}.

To limit the value of the penalty coefficient while controlling
interpenetration, it has been proposed in~\cite{Zavarise-2015} to add
a shift parameter to the penalty contribution, which is incrementally
updated with the residuals obtained when solving the linear system
with an iterative solver. This method results in dynamically shifting
the point from which the penalty potential is applied to the problem
toward the gap-open side, to reduce interpenetration due to contact
without increasing the penalty parameter. This solution has not been
adopted in our work, as it involves correcting the system to be solved
during the solver's iterations, which is not possible when using
black-box solvers such as those from the Hypre library, for
example. On the other hand, a convergence analysis of the penalty
method for 2D and 3D contact problems has been carried out
in~\cite{Chouly-2013}. This study shows that the same mesh convergence
rates with
methods that directly
solution the saddle-point problem arising from the contact formulation are
obtained by choosing $k_{N}$ inversely proportional to $h$ (the mesh size).

\begin{figure}[!h]\centering
\subfloat[][]{ \includegraphics[scale=0.175]{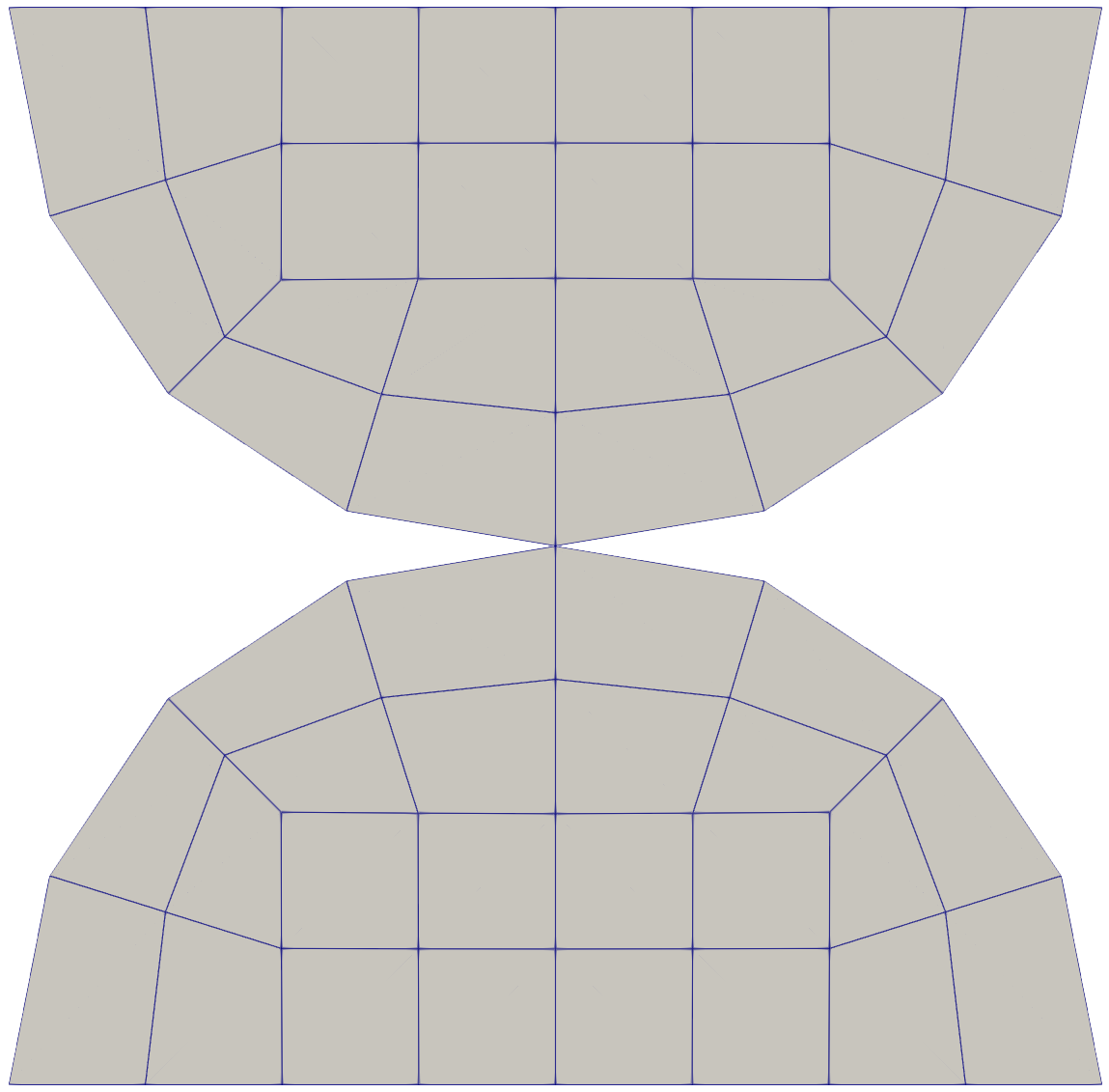} \label{fig:2d_initial_mesh}}
\subfloat[][]{ \includegraphics[scale=0.17]{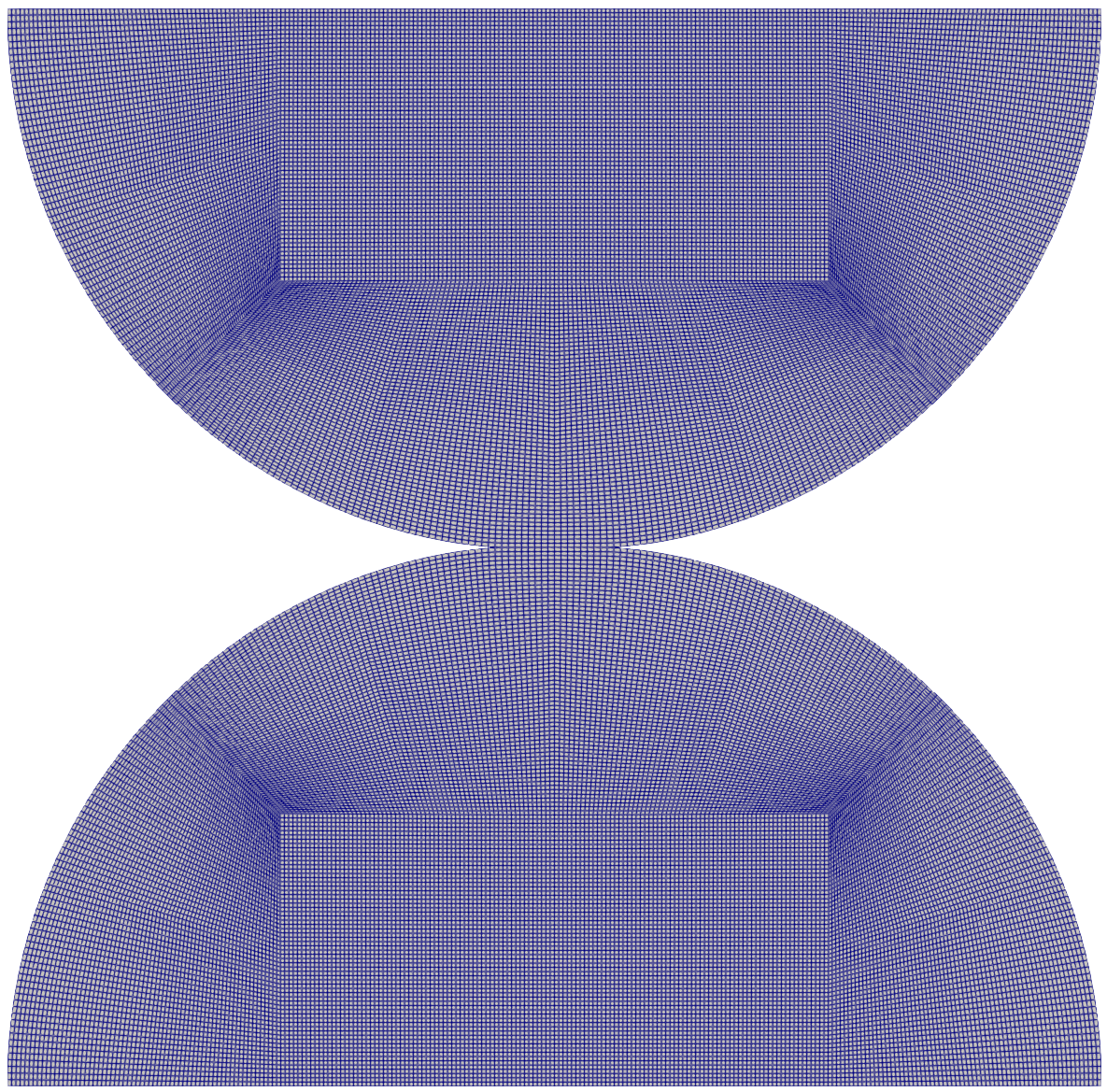} \label{fig:hertz_contact_fig2}}
\caption{2D Hertzian contact problem: initial in (a), and fine in (b)
  quadrilateral meshes}
\label{fig:hertz_contact_meshes}
\end{figure}

Our proposal to choose the penalty coefficient is defined in two
steps. The first step consists in
comparing the results obtained with the penalty method with those
obtained with Lagrange multiplier method. To do this, the contact
forces (see $\textbf{L}^{C}$ in
Problem~\eqref{eqn:generic_contact_problem}) derived from both methods
are compared:
\begin{itemize}
\item $\textbf{L}^{C}(\boldsymbol{\lambda}) = \hat{\textbf{B}}^{T}
\boldsymbol{\lambda}$, for the Lagrange multiplier method, where
$\boldsymbol{\lambda}$ denoted the Lagrange
multiplier solution
\item $\textbf{L}^{C}(\textbf{U}) = k_{N} \hat{\textbf{B}}^{T}
\left( \hat{\textbf{D}} - \hat{\textbf{B}} \textbf{U} \right)$ for the
penalty method
\end{itemize}
Lagrange multiplier contact forces are obtained thanks to simulations
performed with the finite elements software Cast3M~\cite{CEA-2023},
developed by CEA (French Alternative Energies and Atomic Energy
Commission), that is based on a Lagrange multiplier method for contact
mechanics. Penalized contact forces are obtained with our penalized
contact solver implemented in MFEM.  It should be noticed that for all
these comparisons, the input mesh used in Cast3M and MFEM, presented in Figure~\ref{fig:hertz_contact_fig2}, is
identical. Thus, we consider a 2D conforming fine
discretized $Q_{1}$ (bilinear test functions) mesh with a mesh step of
around $h \simeq R / 2^{8} = 2^{-7}$ m, i.e. 200 thousands of nodes in total (100 thousands of nodes in each
half-disk). This mesh is directly
obtained in Gmsh mesher through 6 steps of uniform and isotropic refinement of an initial coarse mesh reported in Figure~\ref{fig:2d_initial_mesh}. Hence, the refined mesh is geometry-fitted. Here, for the 2D meshes, the order of the super-parametric transformation is chosen equal to 10.

\begin{figure}[!h]\centering
\subfloat[][]{ \includegraphics[page=1,clip, trim=0.03cm 0.12cm 0.08cm 0.03cm,width=6.9cm]{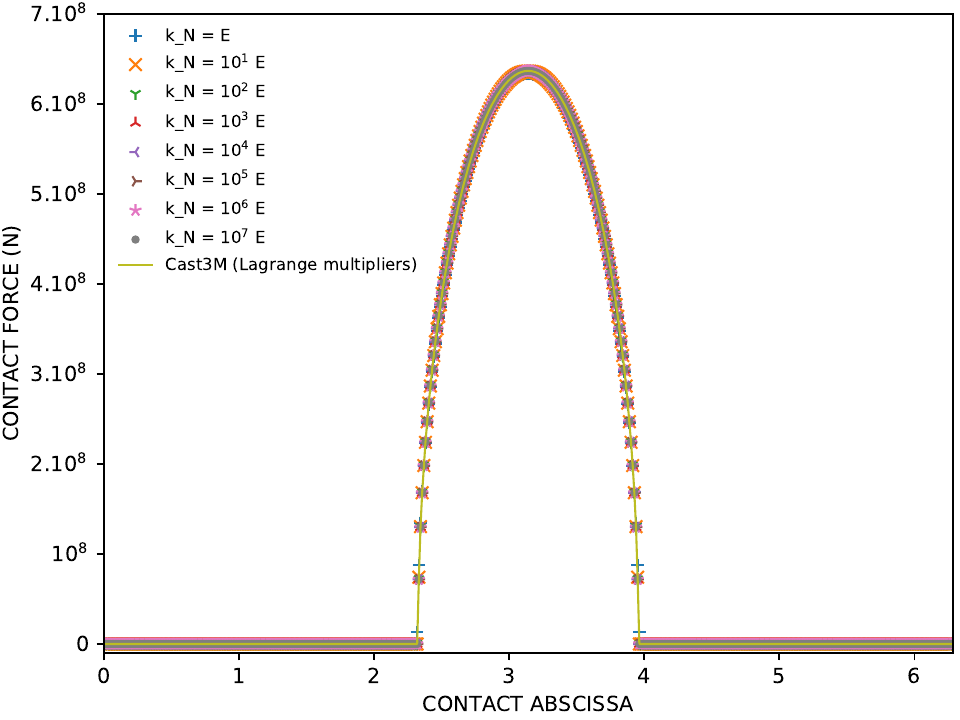} \label{fig:penaltyparam1}}
\newline
\subfloat[][]{ \includegraphics[page=2,clip, trim=0.03cm 0.12cm 0.08cm 0.03cm,width=7cm]{figures/penaltyparam.pdf} \label{fig:penaltyparam2}}
\subfloat[][]{ \includegraphics[page=3,clip, trim=0.03cm 0.12cm 0.08cm 0.03cm,width=7cm]{figures/penaltyparam.pdf} \label{fig:penaltyparam3}}
\caption{Contact force distribution over the potential contact zone
  (half-disks of Hertz). General overview in (a). Zoom on the contact center zone in (b) and on the contact status change area in (c).}
\label{fig:penaltyparam123}
\end{figure}

Figure~\ref{fig:penaltyparam123} plots the contact forces obtained
with our MFEM penalized contact algorithm for different values of the
penalty parameter, from $k_{N} =$E~(N.m$^{-1}$) to $k_{N}= 10^7 \cdot$
E~(N.m$^{-1}$) by steps
of a factor of 10 (without modifying the input mesh), as well as the
contact forces obtained using the Lagrange multiplier method,
considered as the reference solution. Here, $\alpha$ is equal to $0.15$.

As expected, the penalized forces converges towards the Lagrange
multiplier forces when the penalty parameter increases. Zooms made in
Figures~\ref{fig:penaltyparam2} and~\ref{fig:penaltyparam3} clearly
show a poor accuracy for $k_{N} = $~E~(N.m$^{-1}$) and $k_{N} = 10 $
E~(N.m$^{-1}$). Thus,
a penalty coefficient of at least $10^{2}\cdot$E~(N.m$^{-1}$) is required. In
addition, the choice $k_{N} = 1. \dfrac{\text{E}}{h} = 128 \cdot$ E~(N.m$^{-1}$) derived
from~\cite{Chouly-2013}, is perfectly in line with this conclusion and
might be suitable.
\\

\begin{figure}[!h]\centering
\includegraphics[page=1,clip, trim=0.03cm 0.12cm 0.08cm 0.03cm,width=10cm]{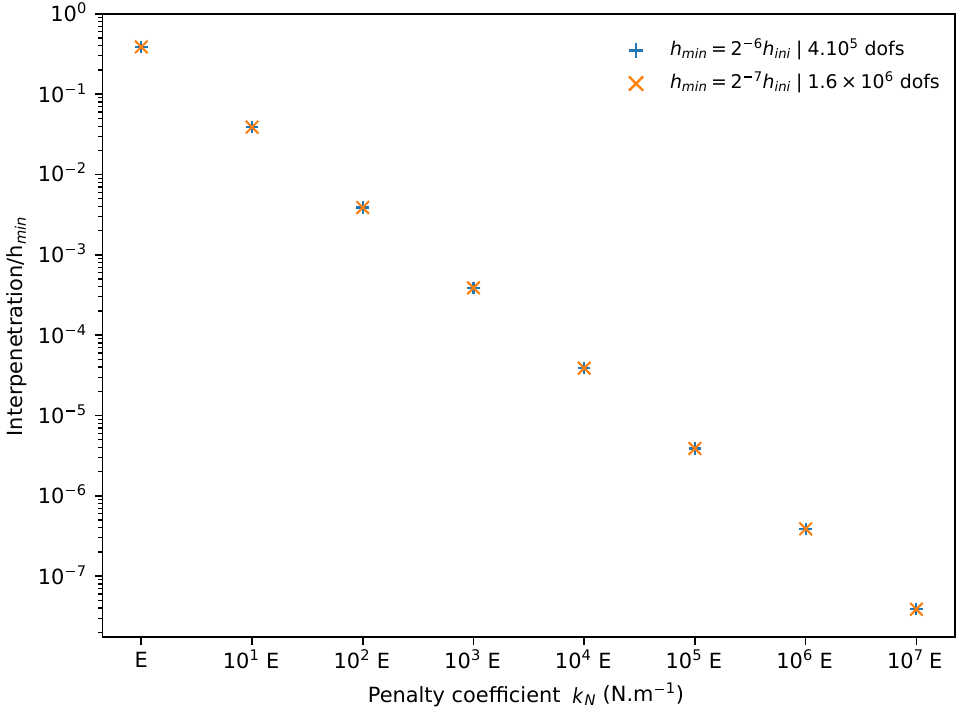}
\caption{2D Hertzian contact problem~: influence of $k_{N}$ and mesh step $h_{min}$ on the interpenetration of the surfaces.}
\label{fig:penaltyparam4}
\end{figure}

In addition, Figure~\ref{fig:penaltyparam4} illustrates the influence of the mesh step and the penalty coefficient on the interpenetration of contact surfaces for two uniformly refined meshes, refined 6 (in blue) and 7 (in orange) times from the same initial mesh (with a mesh step of $h_{ini}$). 
On the one hand, this figure confirms that, for the same finite element simulation, increasing the penalty coefficient logically tends to decrease the interpenetration of the surfaces. On the other hand, this figure allows us to observe that, for the same penalty coefficient, the ratio of interpenetration to mesh step remains constant for the two different meshes.
This suggests that increasing the penalty coefficient is not the only solution to reduce interpenetration. Indeed, decreasing the mesh step leads, for a fixed penalty coefficient, to improve the accuracy of the solution. This result can be justified by the fact that the more nodes there are on the contact boundaries, the lower the nodal value of the contact force on each of these nodes, and the lower the interpenetration resulting from the same penalty. Thus, choosing a constant penalty coefficient during the AMR iterations means reducing the interpenetration with the mesh step. \\

%
%

\begin{figure}[!h]\centering
  \subfloat[][]{\includegraphics[scale=0.2]{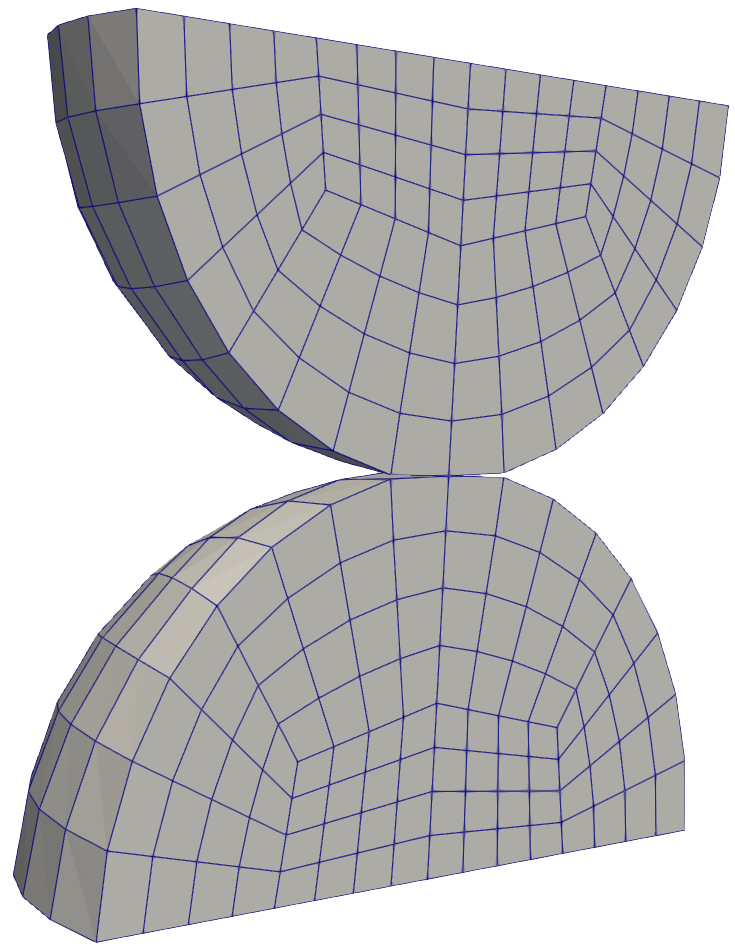} \label{fig:3d_initial_mesh_r2}}
  \hspace*{0.5cm}
  \subfloat[][]{\includegraphics[scale=0.25]{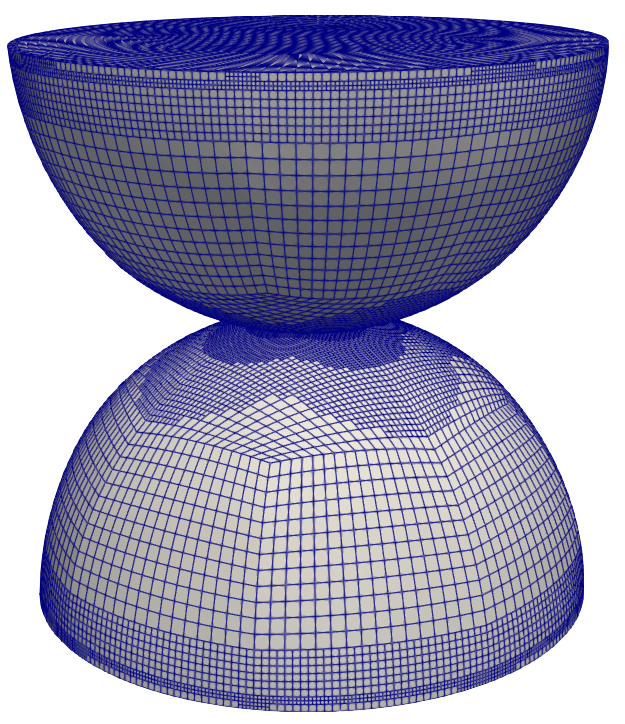} \label{fig:hertz_contact_meshes}}
\caption{3D Hertzian contact problem: (a) cross-sectional view of the
  initial mesh, (b) final mesh for AMR Combination 1 with
$e_{\Omega} = 2 \%$ and $k_N = 10^{2}$ E (N.m$^{-1}$) - 6th-order transformation.}
\end{figure}

The second step consists in comparing strategies for setting $k_{N}$ during a
3D AMR simulation with respect to several quantities of interest:
the number of elements $N_{E}$ and degrees of freedom (DOFs)
$N$,
the total and solver times, $T_{tot}$ and
$T_{solver}$ respectively. Here, AMR Combination 1 with
$e_{\Omega} = 2 \%$ is considered. The initial mesh is built in the Gmsh mesher from the mesh shown in
Figure~\ref{fig:3d_initial_mesh} (in cross-section view), to which a
uniform level of refinement has been applied (by the mesher). The
resulting mesh, presented in Figure~\ref{fig:3d_initial_mesh_r2} (in
cross-section view), contains 2484 nodes and 2048 finite elements and
is fitted to the contact curve boundaries. 
Here, for the 3D meshes, the order of the super-parametric
transformation is chosen equal to 6. The final mesh obtained is shown in
Figure~\ref{fig:hertz_contact_meshes}. Note that here, as for the sequel, $\alpha$ is equal to $0.015$.
Let us underline that the revolution asymmetry of the resulting refined mesh
with respect to the vertical axis is a
direct consequence of the initial mesh that does not respect a
revolution symmetry around this axis.

Simulations were performed on four
calculation nodes of the CCRT/CEA Topaze supercomputer (AMD Milan
processors) each hosting 128 cores per node, so 512 cores in
all. Table~\ref{tab:kn_study_zz_glob_3d}
summarizes some results obtained for the different choices of
$k_N$.
The last choice correspond an adaptation to nodal contact forces with AMR of the suggestion of~\cite{Chouly-2013} where $k_N$ depends on $h_{\min}$, the minimum
discretization size of the current AMR iteration. This results in increasing
$k_N$ during the AMR iterations, which implies even greater
solution accuracy (lower interpenetration) for refined meshes. Let us
underline that keeping $k_N$ constant during AMR iterations already results in
reducing the interpenetration as the mesh step decreases (since the
values of the
nodal forces also decrease).

\begin{table}[!h]
\centering
\begin{tabular}{|C{2.2cm}|C{1.8cm}|C{1.9cm}|C{1.5cm}|C{1.7cm}|}
\hline  $k_{N}$ (N.m$^{-1}$) & $N_{E}$ & $N$ & $T_{tot}$ (s) & $T_{solver}$ (s)  \\
\hline  E & $7.85 \times 10^{5}$ & $1.377 \times 10^{7}$ & $452$ & $311$  \\
\hline  $10$ E & $7.84 \times 10^{5}$ &$1.375 \times 10^{7}$ & $426$ & $308$  \\
\hline  $10^{2}$ E & - & - & $427$ & $309$  \\
\hline  $10^{3}$ E & - & - & $436$ & $304$  \\
\hline  $10^{4}$ E & - & - & $445$ & $305$  \\
\hline  $10^{5}$ E & - & - & $441$ & $306$  \\
\hline  $10^{6}$ E & - & - & $426$ & $307$  \\
\hline  $10^{7}$ E & - & - & $437$ & $303$  \\
\hline  $\dfrac{\text{E}}{h_{\min}}$ & - & - & $431$ & $306$  \\
\hline 
\end{tabular}
\caption{Influence of the penalty parameter $k_N$ on the AMR
  simulation with a target  global error $e_\Omega=2\%$ - Average times over 3 calculations - The '-' symbol indicates that the value is identical to the last entry in the column.}
\label{tab:kn_study_zz_glob_3d}
\end{table}

All simulations ended after 5 refinement
steps, so that the final mesh step
$h_{fin}$ is around $2^{-5} h_{ini}$ with $h_{ini}$ the initial mesh
step. For $k_{N}$ greater than $10 \cdot$ E~(N.m$^{-1}$),
the results in terms of AMR are similar.
In addition, the computation times vary only slightly with the chosen
penalty coefficient. Indeed, variations of around 6 \% are observed
for the total computation time, and not exceeding 3 \% for the solver
time. Thus, the results tend not to depend too much on the penalty
coefficient as soon as this latter is superior or equal to
$10$ E~(N.m$^{-1}$). By choosing a variable penalty parameter as
defined in Table~\ref{tab:kn_study_zz_glob_3d} (last row), it is observed that $k_{N}$ exceeds
$10$ E~(N.m$^{-1}$) starting from the fourth AMR iteration and seems also suitable. \\

From this two studies, it was decided to set the penalty coefficient to be
constant at $k_{N} = 10^{4} \cdot$ E~(N.m$^{-1}$) to achieve a satisfactory compromise
between surface interpenetration and simulation convergence, even for
small elements. Indeed, choosing $k_N = 10^{4} \cdot$ E~(N.m$^{-1}$) does not seem any more
costly in terms of calculation time than $k_N=10^{2} \cdot$ E~(N.m$^{-1}$), and reduces
interpenetration even further.
Increasing the penalty coefficient with the mesh step is not
required here as the penalization directly affects the nodal forces.


\subsection{Comparison with Hertz analytical solution}

To assess the accuracy of the proposed solution strategy, it is relevant to
compare the results obtained with those derived from the Hertz's
theory (see~\cite{Hertz-1882} for further details). To stick to the
scope of Hertz's theory, the following conditions must be satisfied:
\begin{itemize}
\item solids have elastic, homogeneous and isotropic mechanical
  behaviour;
\item the infinitesimal strain theory is applied; 
\item the contact surfaces are small compared to the whole solids and the deformations outside these surfaces are negligible compared to those on them;
\item no contact friction.
\end{itemize}

Thus, for this analysis, the 2D test case
is under consideration and the fixed
Dirichlet boundary conditions ($u_{D} - \dfrac{\delta_{0}}{2} = 0.015~R$) respect the criteria of the infinitesimal strain theory. In addition, we still
rely on the Combination~1 for the AMR process.

\begin{figure}[!h]
\includegraphics[page=1,clip, trim=0.1cm 0.25cm 0.1cm 0.04cm,width=12cm]{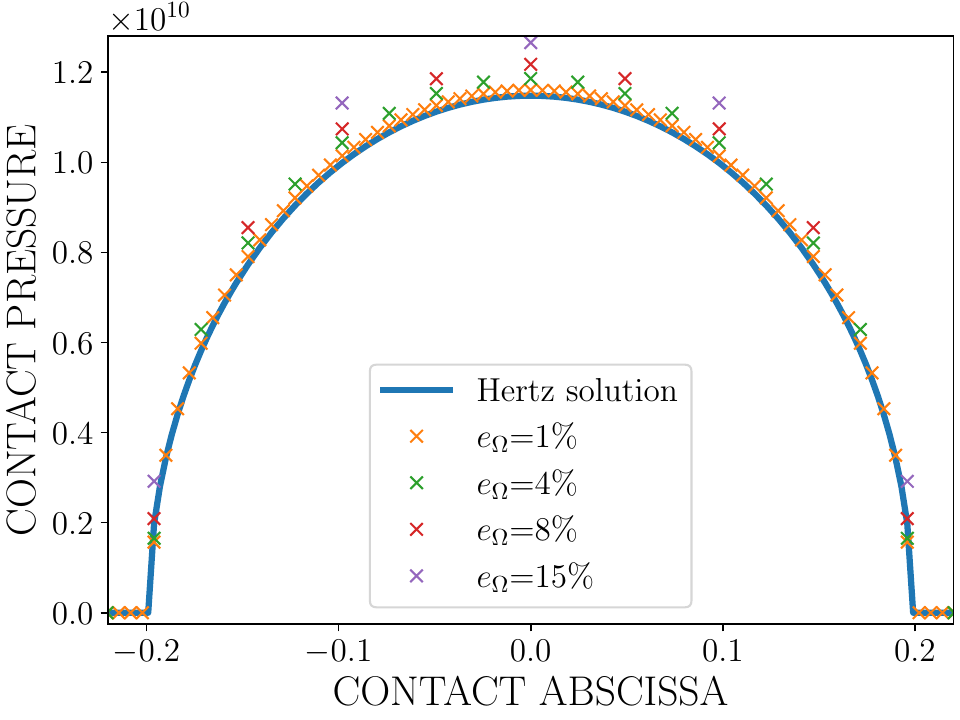}
\caption{Contact pressure distribution over the contact zone
  (half-disks of Hertz).}
\label{fig:hertz_ana_2d_presscontact}
\end{figure}

Hertz's theory defines the following contact pressure distribution: \begin{equation}
p \left( r \right) = p_{\text{O}} \sqrt{a^{2} - r^{2}} / a \; ; \; r \leq a
\end{equation}
where $a$ is the characteristic contact distance, $r$ the (positive)
distance from the center O of the contact zone, and
$p_{\text{O}} = p \left( \text{O} \right) $ the maximum contact
pressure at center O. To calibrate the Hertz solution for our test
case, a highly refined mesh has been used to perform a reference
solution, which yields $a = 0.199$ and $p_{\text{O}} = 11.5$
GPa. These values are then used to define the Hertz pressure
distribution, which is compared in
Figure~\ref{fig:hertz_ana_2d_presscontact} with the contact
pressures obtained via the finite
element solution
$p = - \sigma_{C} = - \left( \boldsymbol{\sigma} \textbf{n}
\right) \cdot \textbf{n}$ on $\Gamma_C$. In this Figure, the origin of the contact
abscissa is chosen as the center O of the contact zone. The results
clearly show that by decreasing the user-specified accuracy $e_{\Omega}$ for
the AMR application, the finite element solution is refined in the
contact area and tends towards the
closed-form expression derived from Hertz theory. Moreover, this
example shows the efficiency of the super-parametric transformation to
capture the curved boundary during the hierarchical AMR iterations.

\subsection{On the AMR strategies} \label{sec:AMR_criteria_analysis}

\begin{figure}[!h]
  \centering
\subfloat[][]{ \includegraphics[page=4,clip, trim=0.12cm 0.25cm 0.15cm
  0.02cm,width=0.75\textwidth]{figures/2d_vs_ref_v2.pdf} \label{fig:cvg-meshstep}}\\
\subfloat[][]{ \includegraphics[page=1,clip, trim=0.12cm 0.15cm 0.09cm
  0.06cm,width=0.75\textwidth]{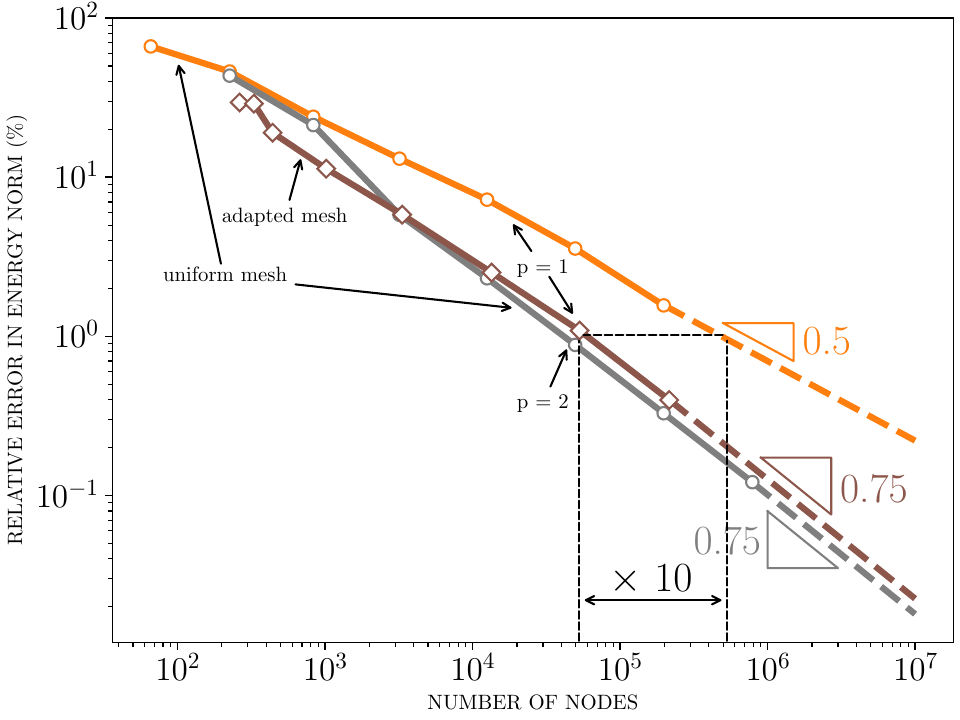} \label{fig:cvg-ddl}}
\caption{Discretization errors (half-disks of Hertz) with respect to the
  mesh size and to
  number of nodes for uniform refinement and AMR approaches.}
\label{fig:2d_vs_ref_v2}
\end{figure}

In order to assess the interest of performing AMR instead of uniform
refinement, convergence in mesh-step and in number of nodes is
performed for the 2D case. To this end, a reference solution (on the
whole surface) is performed on a very fine discretized mesh of linear
finite elements with a mesh step of around $R / 2^{10}$, i.e. 3.2
million nodes in all (1.6 million nodes for each half-disk). 
For this mesh
convergence study, linear or quadratic isoparametric finite elements
are considered.  Indeed, as mentioned in~\cite{Drouet-2015} and the
references cited therein, the contact solution is at best in the
Hilbert space $H^{5/2}$.  Hence, only finite element methods of order
one and of order two are really of interest. To
evaluate the approximation error on the computational domain, uniformly
isotropically pre-refined (by the Gmsh mesher) meshes are used.
An example of
an initial mesh with first-order elements is given in
Figure~\ref{fig:2d_initial_mesh}. 
In its turn, the AMR
method is applied on linear finite elements only. In addition, as mesh
convergence studies involves global errors (on the whole domain), AMR
Combination 1 is considered here. Figure~\ref{fig:2d_vs_ref_v2}
shows the mesh convergence of uniformly refined solutions as well as
the convergence in mesh nodes for uniformly and locally (with AMR)
refined solutions.



On the one hand, it is observed that the uniformly refined solutions
respect the convergence rates predicted by finite element theory in
the presence of contact between two solids, see~\cite{Wohlmuth-2012,
  Drouet-2015} for example. For the energy norm, these convergence rates are
in $\mathcal{O}\left( h^{min\left( k, \tau-1 \right)} \right) $ according to
the mesh size and in
$\mathcal{O}\left( N^{-\min\left( k/m, (\tau-1)/m \right)} \right) $
with respect to
the number of DOFs (directly proportional to the number
of mesh nodes). In these formula, $k$
denotes the degree of interpolation of the basis functions, $\tau$ the
regularity of the solution (here $\tau=5/2$), and $m$ the
dimension of the space.
Hence, for contact problems solved by finite elements of order
greater than 1, sub-optimal convergence rates (here in 3/2) are
obtained, see Figure~\ref{fig:cvg-meshstep}. 

On the other hand, first-order solutions obtained with AMR exhibit
better convergence rates in DOFs than solutions obtained with a uniformly refined
mesh, see Figure~\ref{fig:cvg-ddl}. A superconvergence of the
first-order AMR solutions is observed, as they perform as well as
uniformly refined second-order solutions. As expected, AMR errors are
smaller than for uniformly refined first-order solutions for a given
number of DOFs. Moreover, the number of DOFs is at least approximately
10 times smaller for the AMR solution for a given error producing
sufficiently accurate solutions. This substantial gain justifies the
interest in implementing an AMR strategy.
\\

The ESTIMATE-MARK-REFINE approach considered in our study, combining
quantitative detection and stopping criteria, has recently proven its
efficiency in a solid mechanics context~\cite{Koliesnikova-2021} but
without contact. Thus, the enrichment of the currently implemented
h-adaptive method in MFEM with this strategy is evaluated on the
Hertzian contact problem hereafter. For this case study, the 3D half-spheres are brought into contact and the two AMR combinations introduced above are considered. \\
\textbf{Remarks:} Since a very fine, uniformly refined 3D mesh is very
costly in terms of computation time and memory footprint, error
results are not compared here to a reference solution.
The initial mesh for calculations (performed here with 512 MPI processes) is shown in Figure~\ref{fig:3d_initial_mesh_r2}. \\

\begin{figure}[!h]\centering
\subfloat[][]{ \includegraphics[scale=0.155]{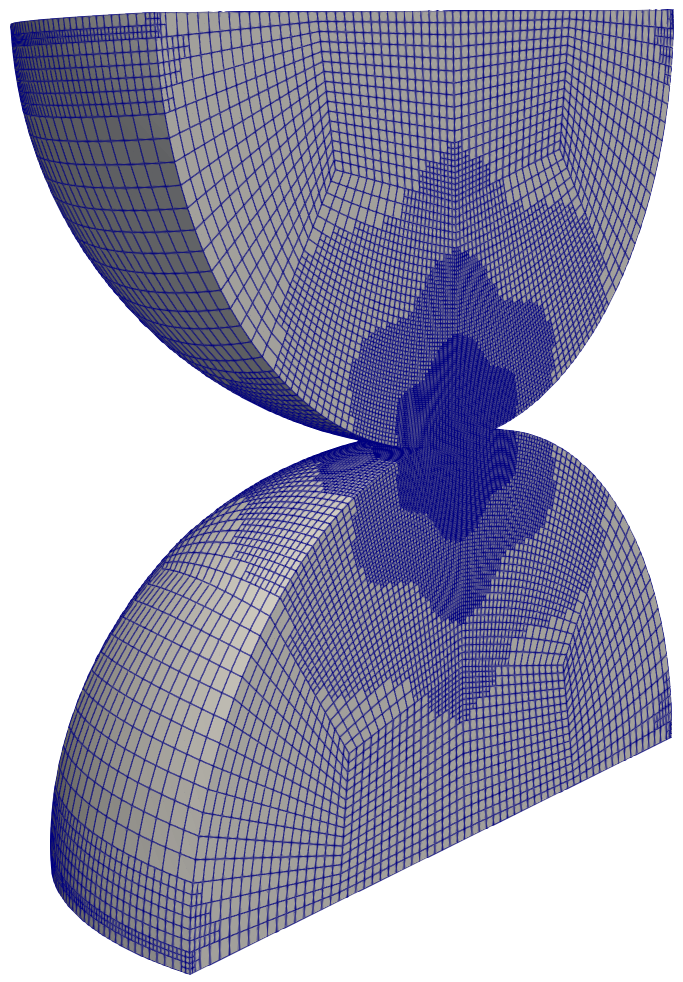} \label{fig:3d_amr_err2_glob_fig1_glob}}
\hspace*{2.0\baselineskip} 
\subfloat[][]{ \includegraphics[scale=0.128]{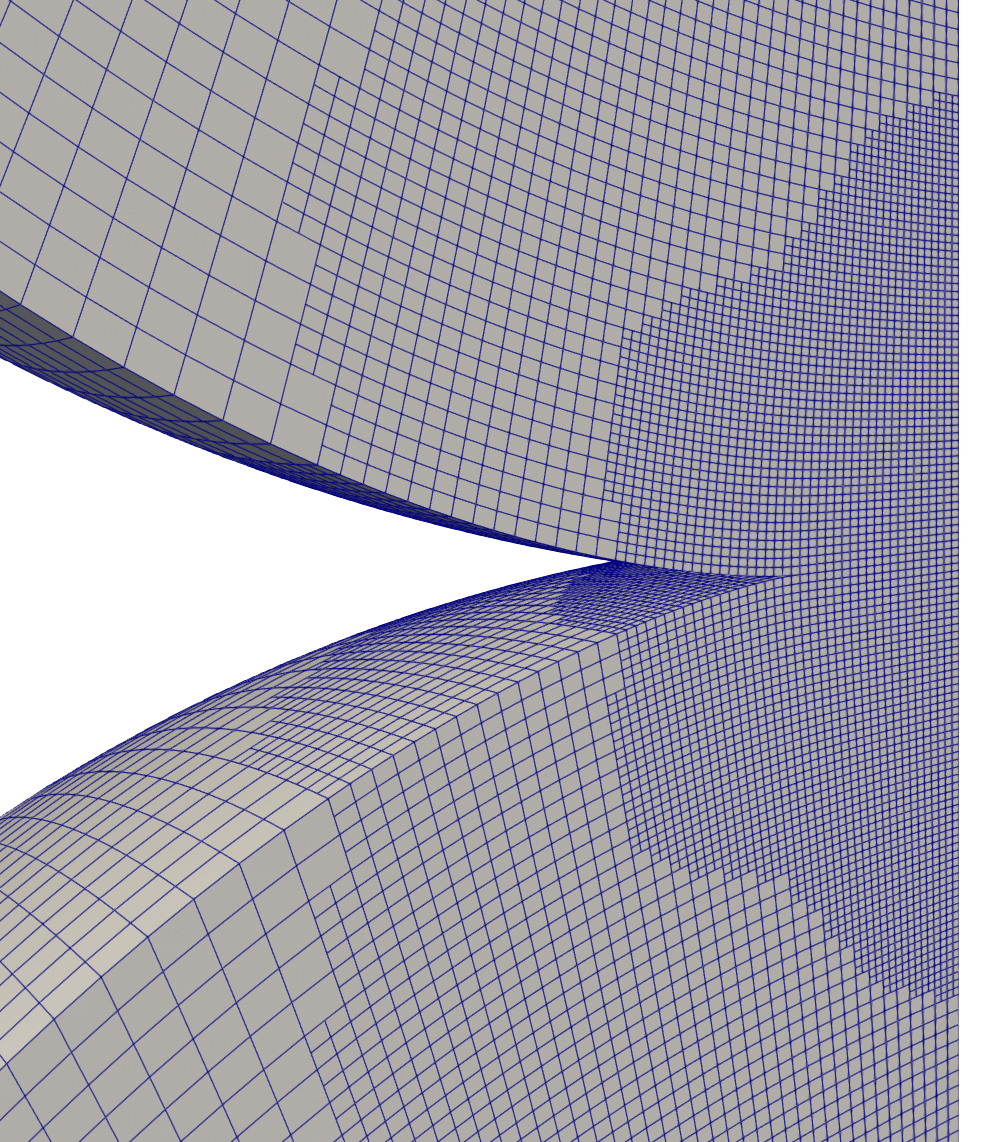} \label{fig:3d_amr_err2_glob_fig1_zoom}}
\caption{3D Hertzian contact - AMR Combination 1 ($e_{\Omega} = 2
  \%$). Cross sectional view in (a) - and zoom around the contact center zone in (b) - 6th-order transformation.}
\label{fig:3d_amr_err2_glob_fig1}
\end{figure}

Table~\ref{tab:zz_glob_3d} details the results obtained with
AMR Combination 1 for different values of $e_{\Omega}$ from $2 \%$ to
$1 \%$ while Figure~\ref{fig:3d_amr_err2_glob_fig1} shows the refined mesh obtained for $e_{\Omega} = 2 \%$.

\begin{table}[!h]
\centering
\begin{tabular}{|C{2cm}|C{1.5cm}|C{1.5cm}|C{2cm}|C{3cm}|}
\hline  & $\gamma_{\Omega} (\%)$ &  $\eta (\%)$ &  $N_{E}$ & Number of refinement steps  \\
\hline  $e_{\Omega} = 2 \%$ & $1.84$ & $15.4$ & $7.84 \times 10^{5}$ & 5  \\
\hline  $e_{\Omega} = 1.7 \%$ & $1.14$ & $29.9$ & $2.41 \times 10^{6}$ & 6  \\
\hline  $e_{\Omega} = 1.5 \%$ & $1.08$ & $33.0$ & $3.07 \times 10^{6}$ & 6  \\
\hline  $e_{\Omega} = 1.2 \%$ & $0.976$ & $21.2$ & $4.59 \times 10^{6}$ & 7  \\
\hline  $e_{\Omega} = 1 \%$ & $0.925$ & $16.1$ & $6.27 \times 10^{6}$ & 7  \\
\hline 
\end{tabular}
\caption{3D Hertzian contact - AMR Combination 1 - results in errors and number of AMR iterations}
\label{tab:zz_glob_3d}
\end{table}


\begin{figure}[!h]\centering
\subfloat[][]{ \includegraphics[scale=0.14]{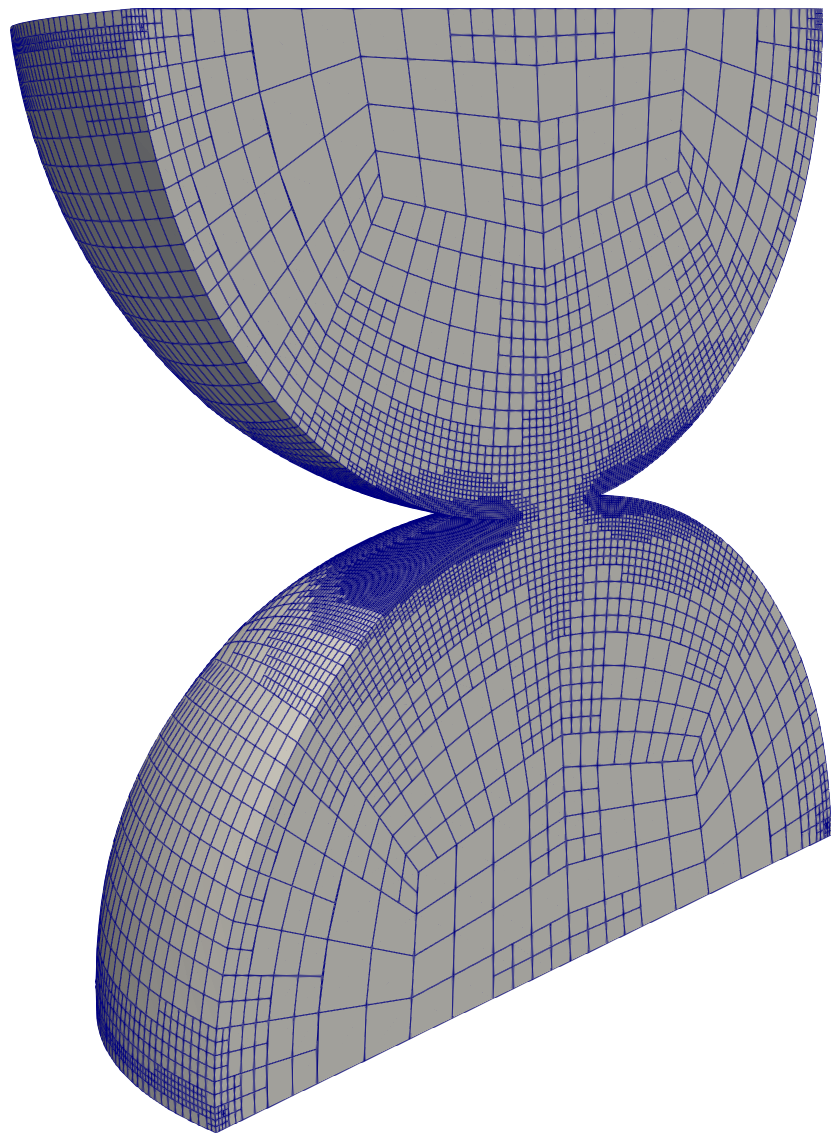} \label{fig:hertzcontact3D_amrloc7_delta0p1_22jui24_mesh_128cpu_glob}}
\hspace*{2.0\baselineskip} 
\subfloat[][]{ \includegraphics[scale=0.129]{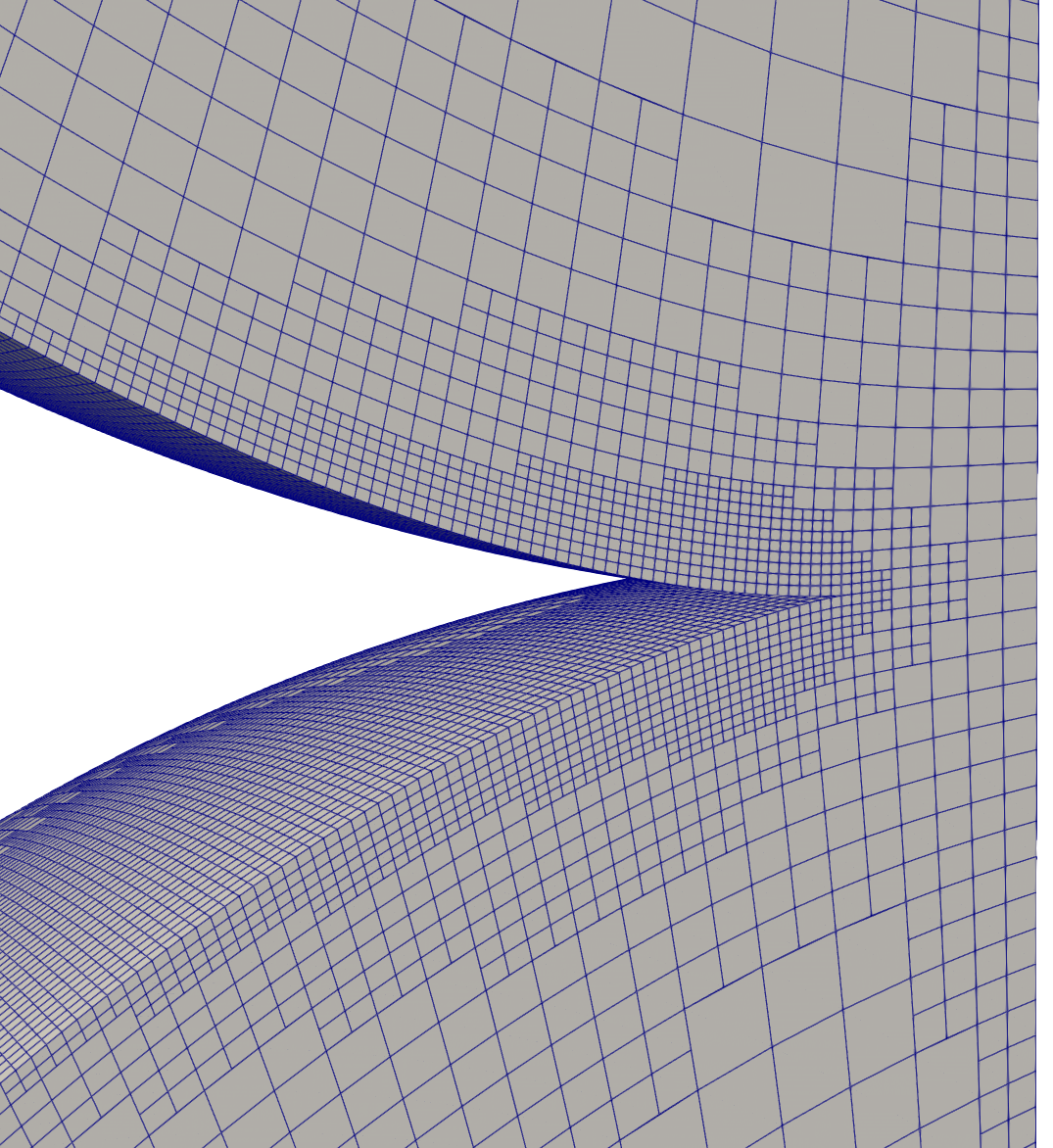} \label{fig:3d_amr_err2_loc_fig1_zoom}}
\caption{3D Hertzian contact - AMR Combination 2 ($e_{\Omega,
    \text{LOC}} = 7 \%$, $\delta = 0,1\%$). Cross sectional view in (a)
  and zoom on the contact center zone in (b) - 6th-order transformation.}
\label{fig:hertzcontact3D_amrloc7_delta0p1_22jui24_mesh_128cpu}
\end{figure}

Similarly, Table~\ref{tab:loc_loc_3d} provides the results obtained
with AMR Combination 2 for $e_{\Omega , \text{LOC}} = 7 \%$ to $2\%$ and
$\delta = 0.1 \%$, while
Figure~\ref{fig:hertzcontact3D_amrloc7_delta0p1_22jui24_mesh_128cpu}
displays the refined mesh obtained for $e_{\Omega , \text{LOC}} = 7 \%$ and $\delta = 0.1 \%$. Here, $e_{\Omega , \text{LOC}}$ is globally chosen to be greater than $e_{\Omega}$ because AMR Combination 2 tends to locally more refine than AMR Combination 1, see Tables~\ref{tab:zz_glob_3d} and~\ref{tab:loc_loc_3d} for the error value of $2\%$ for example.

\begin{table}[!h]
\centering
\begin{tabular}{|C{4.5cm}|C{1.5cm}|C{1.5cm}|C{2cm}|C{3cm}|}
\hline  & $\gamma_{\Omega} (\%)$ &  $\eta (\%)$ &  $N_{E}$ & Number of refinement steps  \\
\hline  $e_{\Omega , \text{LOC}} = 7 \%$, $\delta = 0.1 \%$ & $5.11$ & $0.0537$ & $2.88 \times 10^{5}$ & 5  \\
\hline  $e_{\Omega , \text{LOC}} = 5 \%$, $\delta = 0.1 \%$ & $3.60$ & $0.0182$ & $1.05 \times 10^{6}$ & 6  \\
\hline  $e_{\Omega , \text{LOC}} = 4 \%$, $\delta = 0.1 \%$ & $2.92$ & $0.0330$ & $1.73 \times 10^{6}$ & 6  \\
\hline  $e_{\Omega , \text{LOC}} = 3 \%$, $\delta = 0.1 \%$ & $2.21$ & $0.0108$ & $4.64 \times 10^{6}$ & 7  \\
\hline  $e_{\Omega , \text{LOC}} = 2 \%$, $\delta = 0.1 \%$ & $1.47$ & $0.0363$ & $1.34 \times 10^{7}$ & 7  \\
\hline 
\end{tabular}
\caption{3D Hertzian contact - AMR Combination 2 - results in errors and number of AMR iterations}
\label{tab:loc_loc_3d}
\end{table}

\begin{figure}[!h]\centering
\subfloat[][]{ \includegraphics[scale=0.173]{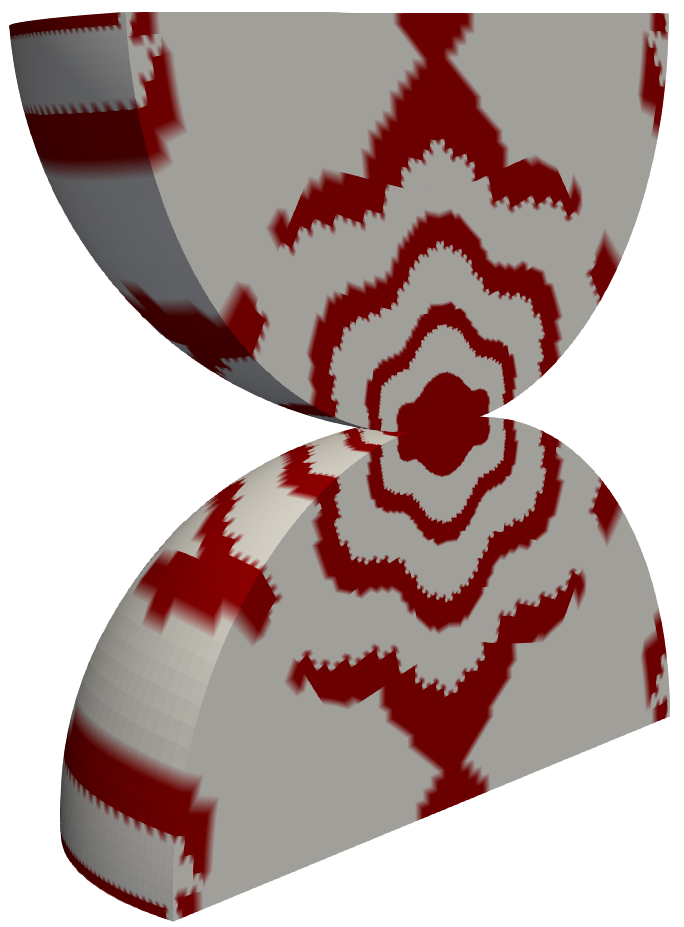} \label{fig:hertzcontact3D_amrglob2_22jui24_markedelements_128cpu_glob}}
\hspace*{2.0\baselineskip} 
\subfloat[][]{ \includegraphics[scale=0.14]{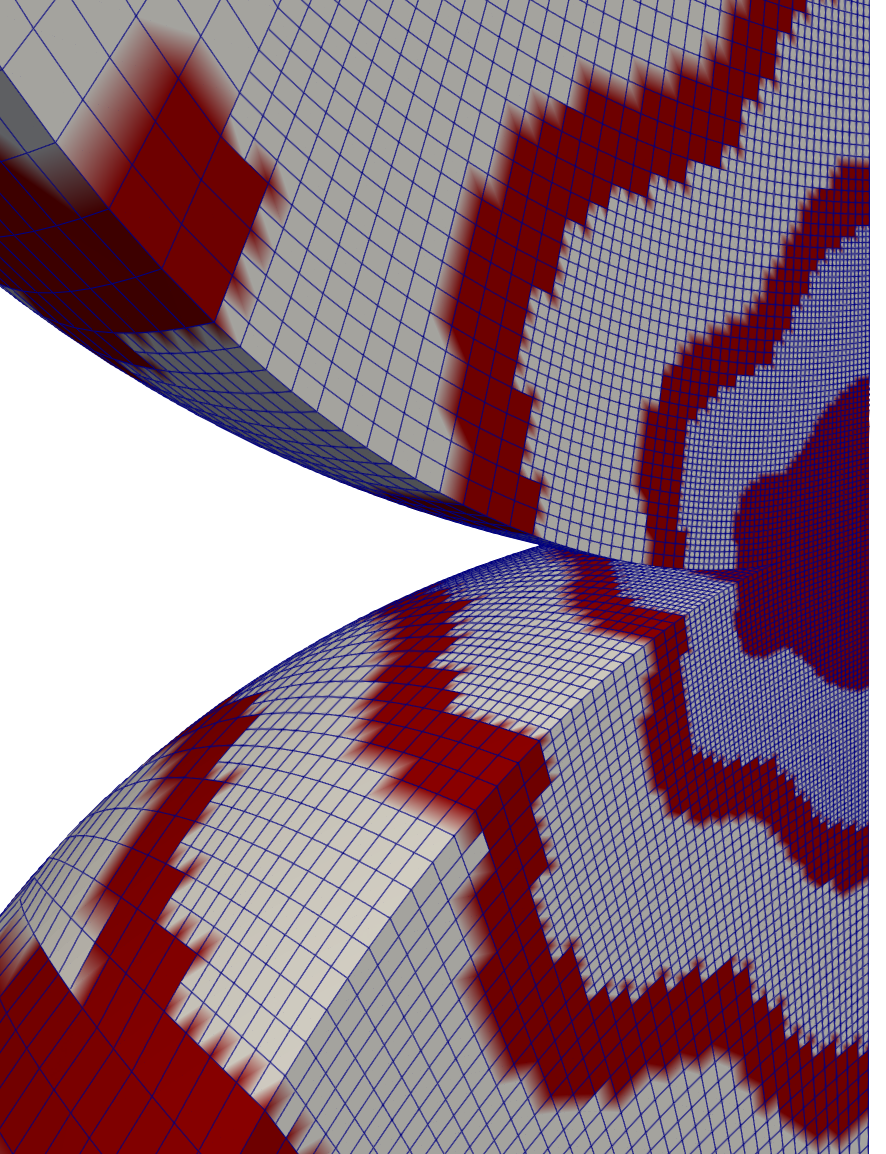} \label{fig:hertzcontact3D_amrglob2_22jui24_markedelements_128cpu_zoom}}
\caption{3D Hertzian contact - AMR Combination 1 ($e_{\Omega} = 2 \%$). Zones still marked for refinement at the end of the AMR process in red. Cross sectional view in (a). Zoom on the contact center zone in (b) - 6th-order transformation.}
\label{fig:hertzcontact3D_amrglob2_22jui24_markedelements_128cpu}
\end{figure}

\begin{figure}[!h]\centering
\subfloat[][]{ \includegraphics[scale=0.14]{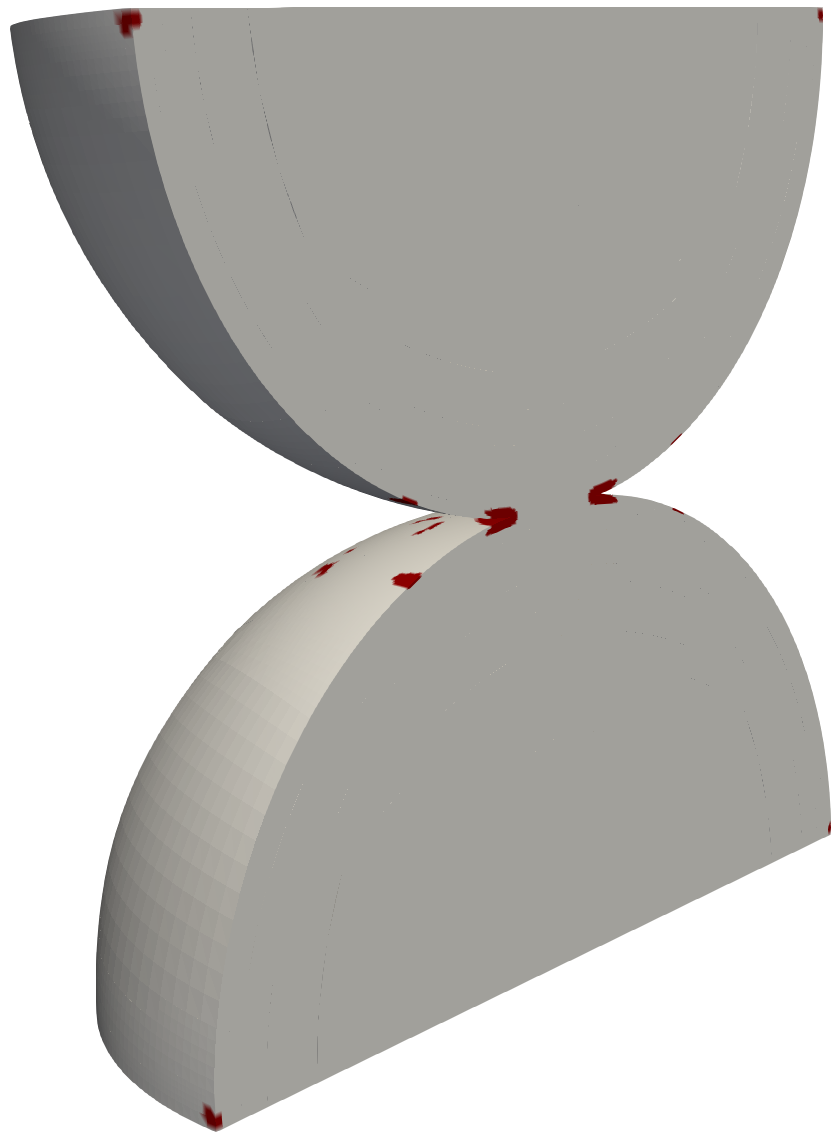} \label{fig:hertzcontact3D_amrloc7_delta0p1_22jui24_markedelements_128cpu_glob}}
\hspace*{2.0\baselineskip} 
\subfloat[][]{ \includegraphics[scale=0.149]{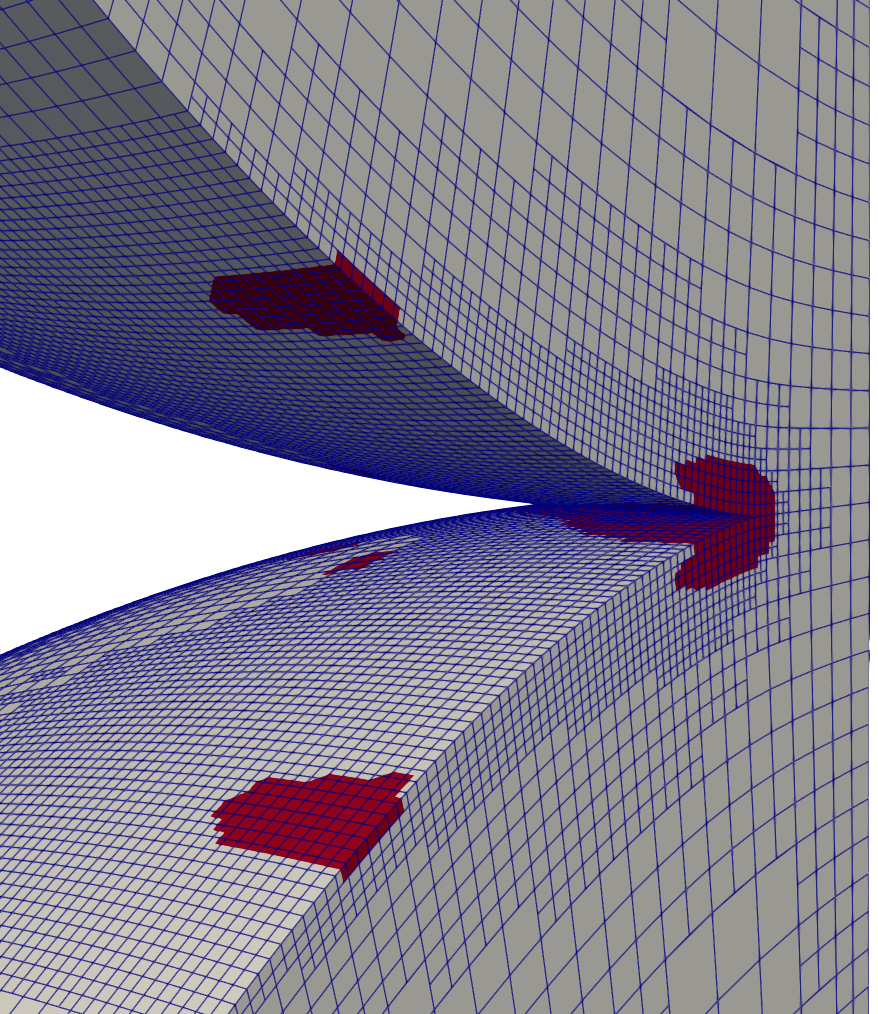} \label{fig:3d_amr_err7_delta_0point1_loc_fig2_zoom}}
\caption{3D Hertzian contact - AMR Combination 2 ($e_{\Omega, \text{LOC}} = 7 \%$, $\delta = 0,1\%$). Zones still marked for refinement at the end of the AMR process in red. Zoom on the contact center zone in (b) - 6th-order transformation.}
\label{fig:hertzcontact3D_amrloc7_delta0p1_22jui24_markedelements_128cpu}
\end{figure}

First of all, we observe
that for AMR Combination 1 (see Figure~\ref{fig:3d_amr_err2_glob_fig1}), mesh refinement is spread in the contact
area, whereas for AMR Combination 2,
Figure~\ref{fig:hertzcontact3D_amrloc7_delta0p1_22jui24_mesh_128cpu} shows that the finest edges
(corresponding to the highest level of refinement) are located in the
zone where the contact status (contact/no contact) changes, which is in
coherence with the results obtained in~\cite{Liu-2017} with the same
LOC optimality criterion. Both tables show that, whatever
the selected combination, the user-prescribed error threshold
$e_{\Omega}$ (resp. $e_{\Omega,
    \text{LOC}}$) is respected, i.e. $\gamma_\Omega < e_\Omega$ (resp. $\gamma_\Omega < e_{\Omega,
    \text{LOC}}$), especially for AMR Combination 2 where this condition is not imposed. In addition, for an equivalent relative
global error $\gamma_{\Omega}$, the mesh obtained contains fewer
elements when using AMR Combination 1, but the use of this combination
does not effectively control the zone where the local prescribed error is not
respected, see $\eta$ in row 1 of Table~\ref{tab:zz_glob_3d} versus in
row 5 of
Table~\ref{tab:loc_loc_3d}. On the contrary, choosing AMR Combination 2 to
apply AMR enables us to more efficiently characterize the critical
regions, see
Figure~\ref{fig:hertzcontact3D_amrloc7_delta0p1_22jui24_markedelements_128cpu}
versus
Figure~\ref{fig:hertzcontact3D_amrglob2_22jui24_markedelements_128cpu},
where the number of AMR iterations is the same. The

Thus, the choice of the relevant combination of AMR criteria depends
on what is of interest to control. On the one hand, if the main aim is
to control the global error only, then choosing Combination 1 enables
to obtain an acceptable mesh with less DOFs. On the other hand, if our priority is to locally handle the
error, then choosing Combination 2 is required.

The Figures~\ref{fig:3d_amr_err2_glob_fig1}
to~\ref{fig:hertzcontact3D_amrloc7_delta0p1_22jui24_markedelements_128cpu}
show that the symmetrical nature of the problem
(identical geometries,
same materials) results in perfectly identical refinement of the two
solids, where the nodes are strictly facing each other.
However, as previously mentioned the initial mesh induces a revolution
asymmetry around the vertical axis of all the refined meshes.
\\

\begin{figure}[!h]\centering
\subfloat[][]{ \includegraphics[scale=0.20]{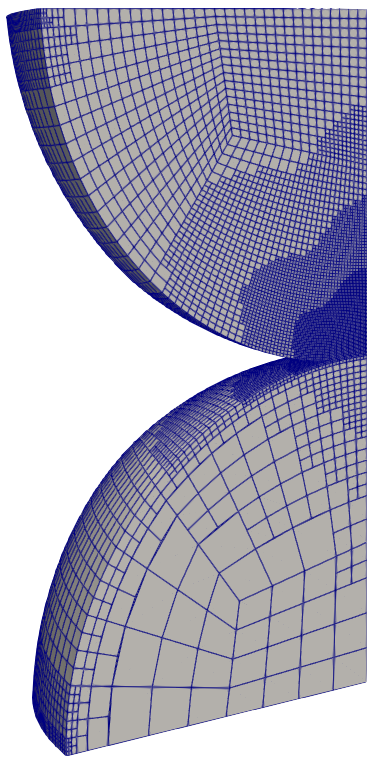} \label{fig:3d_amr_err2_glob_2mat_fig1_glob}}
\hspace*{2.0\baselineskip} 
\subfloat[][]{ \includegraphics[scale=0.189]{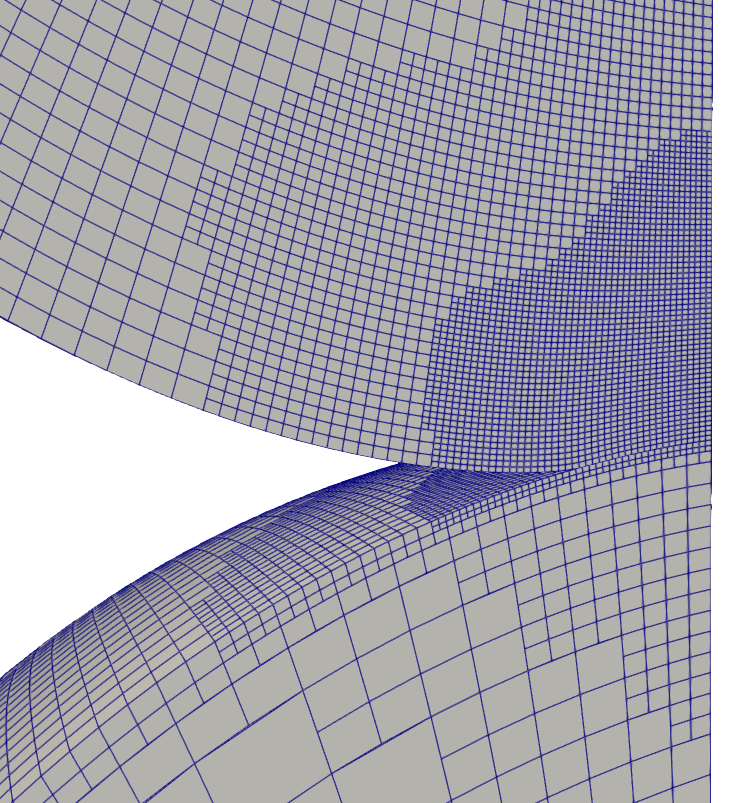} \label{fig:3d_amr_err2_glob_2mat_fig1_zoom}}
\caption{3D Hertzian contact between different materials -
  AMR Combination 1 ($e_{\Omega} = 2 \%$). Cross sectional view in (a) and
  zoom on the contact center zone in (b) - 6th-order transformation.}
\label{fig:3d_amr_err2_glob_2mat_fig1}
\end{figure}

We are also interested in verifying the effectiveness of the AMR
strategy employed in the case of contact between two solids with
different mechanical properties. Simulations were performed with two
solids presenting a contrast of $10^{4}$ in Young's modulus (the Poisson
ratio remains the same). The lower solid
is then more rigid than the upper one, which keeps the material
characteristic previously considered ($E_1= 210$ GPa). The same Dirichlet
boundary conditions than those introduced at the begining of the section were applied
and AMR Combination 1 was chosen. An example of the results obtained
is reported in Figure~\ref{fig:3d_amr_err2_glob_2mat_fig1}. Several
observations are worth noting. Firstly, the AMR process is applied
distinctly between solids in contact, and it can be seen that the most
deformable solid is the most refined overall. Moreover, the mesh
solution obtained shows the effect of the
detection enlargement
criterion~\eqref{eqn:critere_elargi_detection_ensemble_elements_amr},
as all potential contact paired elements are refined in
the same way.

\subsection{Parallel scalability}

To demonstrate the parallel scalability of our method, we have carried
out numerical studies with a varying number of cores for
representative simulations of the 3D half-spheres Hertzian contact problem.

In addition to the total calculation time, the times of the following
main program functions are monitored:
\begin{itemize}
\item assembling the finite element stiffness matrix (\texttt{mfem::ParBilinearForm::Assemble} MFEM function);
\item forming the linear system (\texttt{mfem::ParBilinearForm::FormLinearSystem} MFEM function, equation~\eqref{eqn:assemblage_conforme});
\item solving the matrix system (apply the solver from \texttt{mfem::PenaltyPCGSolver} MFEM class, equation~\eqref{eqn:mechanical_equilibrium_matrix_system}); 
\item refining elements (ESTIMATE-MARK-REFINE AMR approach, \texttt{mfem::MeshOperator::Apply} MFEM function); 
\item load balancing the whole mesh after refinement (including the \texttt{mfem::ParNCMesh::Rebalance} MFEM function).
\end{itemize}

First of all, AMR Combination 1 is applied to the benchmark. The
threshold $e_{\Omega}$ is set to $2 \%$ and the resulting mesh has
been already shown in the previous section, see Figure~\ref{fig:3d_amr_err2_glob_fig1}. This mesh is obtained
after 5 refinement steps and the discrete solution $\textbf{U}$ to be solved 
is of size $2.3 \times 10^{6}$ (number of conforming DOFs).

\begin{figure}[!h]
\centering
\includegraphics[page=1,clip, trim=0.1cm 0.12cm 0.1cm 0.1cm,width=0.8\textwidth]{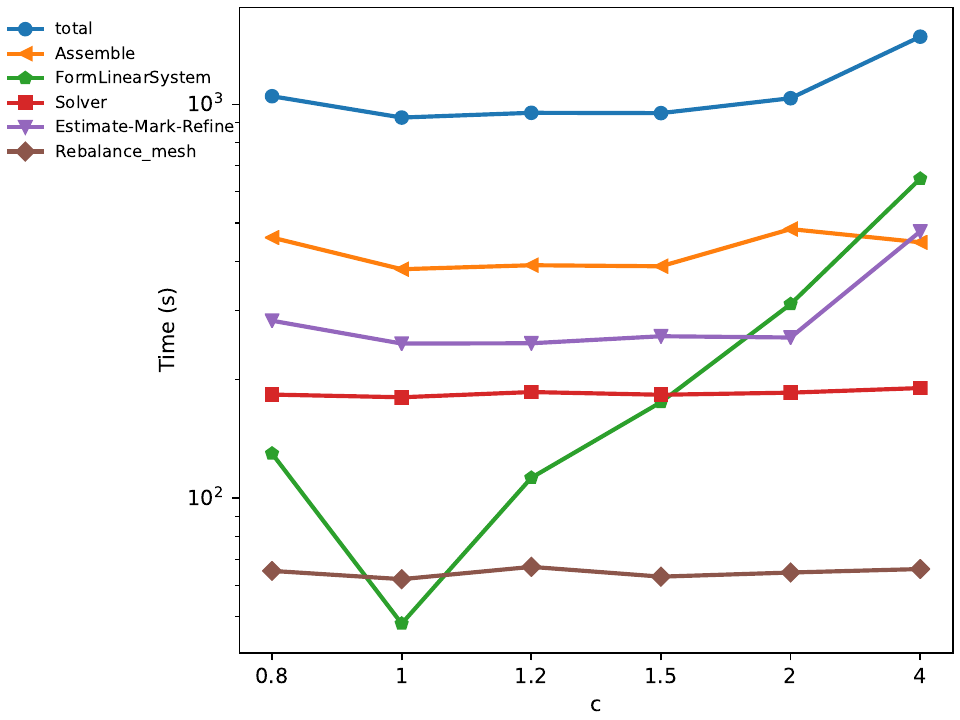}
\caption{Computation times, scan for different values of $c$ - 3D Hertzian contact - AMR Combination 1 ($e_{\Omega} = 2 \%$).}
\label{fig:allKernel_compcoeff_512}
\end{figure}

To set the coefficient $c$ in the load balancing
formula~\eqref{eqn:nb_contact_procs}, calculations are performed by
varying this coefficient for a constant number of MPI tasks. The
results obtained are presented in
Figure~\ref{fig:allKernel_compcoeff_512} for 512 MPI tasks, and show that setting this
parameter to 1 seems to be an optimal choice. Furthermore, since in our strategy MFEM
distributes calculations according to the specified partition of
elements, these results demonstrate that balanced partitioning of the
mesh is essential in optimizing the algorithm's performance.
Other calculations not presented
here have shown that choosing a different number of MPI processes
produces time measurements leading to the same conclusion for the
choice of $c$.

\begin{figure}[!h]\centering
\subfloat{ \includegraphics[page=1,clip, trim=0.7cm 0.15cm 1.5cm 1.35cm,width=7cm]{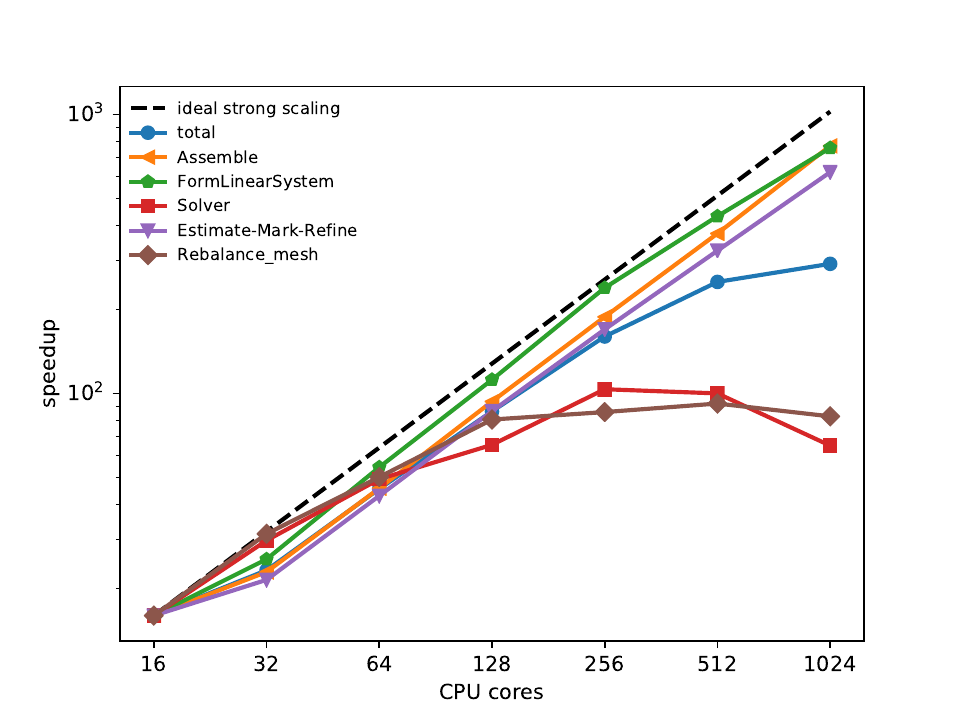}}
\hspace*{0.1\baselineskip} 
\subfloat{ \includegraphics[page=3,clip, trim=0.65cm 0.15cm 1.5cm 1.35cm,width=7cm]{figures/results_amrglob1p5_21mar24_4coeurs_par_proc_HypreILUHyprePCG}}
\caption{Performance results of the contact-AMR-HPC algorithm and its
  main steps - 3D Hertzian contact - AMR Combination 1 ($e_{\Omega} =
  1.5 \%$) - HyprePCG solver with
  HypreILU preconditioner - number of final mesh conforming DOFs: $9.0 \times 10^{6}$.}
\label{fig:allKernelCompPourcentage_amrglob2}
\end{figure}

\begin{figure}[!h]
\centering
\includegraphics[page=3,clip, trim=0.65cm 0.15cm 1.5cm 1.35cm,width=0.7\textwidth]{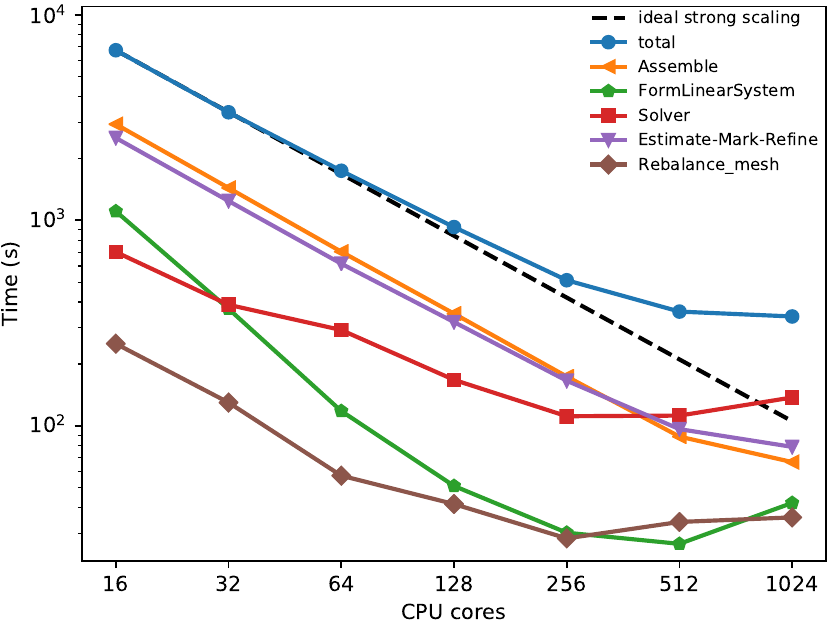}
\caption{Performance results of the contact-AMR-HPC algorithm and its
  main steps - 3D Hertzian contact - AMR Combination 2 ($e_{\Omega,
    \text{LOC}} = 4 \%$, $\delta = 0,1\%$) - HyprePCG solver with
  HypreILU preconditioner - number of final mesh conforming DOFs: $5.1 \times 10^{6}$.}
\label{fig:allKernelCompPourcentage_amrloc2}
\end{figure}

It is of interest to study the parallel performance of our strategy on
larger meshes. This is why, in the scalability analysis reported in
Figures~\ref{fig:allKernelCompPourcentage_amrglob2}
and~\ref{fig:allKernelCompPourcentage_amrloc2}, the contact problem is
solved by applying AMR Combination 1 with $e_{\Omega} = 1.5\%$, and
AMR Combination 2 with $e_{\Omega, \text{LOC}} = 4\%$ and
$\delta = 0.1\%$. The resulting solutions are obtained after 6
refinement steps and the discrete solution  $\textbf{U}$ of the
last 'conforming'
problem is of size $9.0 \times 10^{6}$ and
$5.1 \times 10^{6}$ respectively. Calculations were performed on the
calculation nodes of the CCRT/CEA Topaze supercomputer (AMD Milan
processors) each hosting 128 cores per node. The initial uniformly
refined mesh is designed so that there are always more elements than
the maximal number of
parallel tasks for the study.
The following results are obtained solving the conforming
linear system using Hypre's conjugate gradient solver (HyprePCG),
preconditioned by the incomplete LU factorization strategy available
in Hypre (HypreILU).

During the simulations, it has been observed that the memory footprint
of the calculations was high (mainly due to the distributed mesh
management within MFEM
and constitutes our main limiting
factor. Thus, in order to carry out our strong scaling studies, it
has been decided to allocate 4 cores to each processor
used. Table~\ref{tab:memory_footprint_combinaison1} reports the memory
footprint used for AMR Combination 1 with $e_{\Omega}
= 1.5 \%$ for 256, 512 and 1024 cores. This memory footprint is listed
for calculations using Hypre's ILU preconditioner but also BoomerAMG preconditioner as
well as without any preconditioner. For each simulation, the
percentage of the total memory footprint used is given for the
following steps:
\begin{itemize}
\item the last refinement step ; 
\item the last load rebalancing step ;
\item the steps for solving the matrix system during the last AMR iteration. 
\end{itemize}

\begin{table}[!h]
\centering
\begin{tabular}{|C{3.56cm}|C{0.75cm}|C{0.75cm}|C{0.75cm}|C{1.545cm}|C{1.545cm}|C{0.75cm}|C{0.75cm}|C{0.75cm}|}
\hline Preconditioner  & \multicolumn{3}{c|}{$\emptyset$} &  \multicolumn{2}{c|}{HypreBoomerAMG} & \multicolumn{3}{c|}{HypreILU}  \\
\hline  CPU cores & $256$ & $512$ & $1024$ & $256$ & $512$ & $256$ & $512$ & $1024$  \\
\hline  Total memory footprint (GB) & $662$ & $973$ & $1545$ & $875$ & $1218$ & $671$ & $977$ & $1546$  \\
\hline  Last refinement step ($\%$) & $40.4$ & $34.5$ & $28.2$ & $22.4$ & $21.9$ & $40.7$ & $34.9$ & $28.3$  \\
\hline  Last rebalance step ($\%$) & $8.2$ & $11.0$ & $18.8$ & $6.3$ & $8.6$ & $8.1$ & $10.9$ & $18.8$  \\
\hline  Last AMR iteration - Solver steps ($\%$) & $0$ & $0$ & $0$ & $20.5$ & $17.0$ & $0$ & $0$ & $0$  \\
\hline 
\end{tabular}
\caption{Memory footprint - 3D Hertzian contact - AMR Combination 1
  ($e_{\Omega} = 1.5 \%$)- HyprePCG solver with different
  preconditioner - number of final mesh conforming DOFs: $9.0 \times 10^{6}$ - 6th-order transformation}
\label{tab:memory_footprint_combinaison1}
\end{table}

Note that the simulation for 1024 cores with HypreBoomerAMG
preconditioner did not converge, for a reason that will be explained
later. The increase in memory consumption as a function of the number of
cores is due to the duplication of data structures required for the
ghost layers. Table~\ref{tab:memory_footprint_combinaison1} shows
that, as expected due to the increasing number of DOFs,
the last refinement step consumes a lot of memory space. The last load
rebalancing step is also a major memory space consumer, but less than the last refinement step. Finally, the last solutions of the matrix system induce a
large cost in terms of memory space when the system is solved with a
HypreBoomerAMG preconditioning but not at all with HypreILU or without
preconditioner. This cost in memory space for the solver step with
HypreBoomerAMG preconditioning is due to the storage of the
preconditioning matrix, which is much larger than for HypreILU
preconditioning. Indeed, during the solver steps of the last AMR
iteration, the number of non-zero coefficients in the preconditioning
matrix is about 2 times higher with HypreBoomerAMG than with HypreILU as preconditioner for the calculation using 256 cores, and is 4 times higher for the calculation using 512 cores. When measuring memory, it was observed with HypreBoomerAMG that the memory footprint is larger than the memory footprint due to the theoretical storage of the preconditioning matrix resulting from the number of non-zero coefficients found in Hypre. This could be due to the parallelism of the preconditioning matrix.

Finally, for large 3D calculations and with the aim of reducing computational costs (computation time and memory footprint), a solution could be to reduce the order of the super-parametric transformation or even to implement a the super-parametric transformation solely on the contact boundaries. Nevertheless, this is appropriate for boundaries with curvature that does not vary significantly. \\
On the basis of preliminary calculations, we found for the 3D Hertzian test case that switching from a 6th order to a 3rd order (volume) transformation not only reduces the memory footprint allocated during calculation by a factor of 2, but also divides calculation times by a factor of 5.


For AMR Combination 1,
Figure~\ref{fig:allKernelCompPourcentage_amrglob2} displays the
speedup
in strong scaling, the total time and the individual algorithm steps
times.

Comparing the results obtained with the ideal strong scaling result,
where parallelization would be perfect and would not generate any
overhead, we can observe that the scalability of the algorithm
obtained is satisfactory, particularly up to 512 cores. On the other
hand, a loss of performance is observed regarding the linear system
solving step. The problem has been investigated with the Scalasca
profiling tool~\cite{Bohme-2012}. The loss of scalability is due to
unbalanced communication schemes and large communication volumes
caused by the mesh partitioning.  The reason for this is that the set of
elements owned by each process is composed of elements scattered
throughout the mesh. This leads to large borders between parallel
regions, and therefore extra communication costs. An in-depth
profiling study has shown that the volume of communication is twice as
high as calculations volume.  The most costly operations
are the \texttt{MPI\_Allreduce()} and \texttt{MPI\_Waitall()} reduction and
communication operations, which take up almost $70 \%$ of the time spent
in the solver. So, while the contact-aware partitioning reduces communication 
during contact evaluation, it leads to fragmented MPI domains with large interfacial 
regions, contributing to increased communication overhead during the solver phase.
This topic is beyond the scope of this paper and will be addressed in future work, 
but we can shed some light on the subjects 
we're planning to investigate to reduce this overhead.
An avenue for future research is to divide the finite element mesh into 
larger contiguous parts within each MPI process. This issue will be addressed with 
geometric partitioning techniques to split the mesh elements according to their 
spatial distribution. This kind of partitioning is generally quite cheap 
in term of execution time compared to graph-based approches.  
In this regard, we are currently working on a solution that uses the Space Filling Curve
(SFC) algorithm to generate a mesh partition and create large, parallel regions composed of 
closely grouped elements while ensuring contact paired elements to be on the same process.

A comparative study of solver performances regarding the choice of preconditioner is
presented in Table~\ref{tab:times_combinaison1}. This table shows the
number of iterations and execution times spent in the linear system
solution step for two peculiar solving steps: the last (contact problem) solving step for the
last two AMR iterations. These results are reported for 256, 512 and
1024 CPU cores for HypreBoomerAMG and HypreILU preconditioning, as
well as without preconditioner.

\begin{table}[!h]
\centering
\begin{tabular}{|C{2.8cm}|C{1.6cm}|C{0.8cm}|C{0.8cm}|C{0.8cm}|C{0.89cm}|C{0.89cm}|C{0.89cm}|C{0.8cm}|C{0.8cm}|C{0.8cm}|}
\hline \multicolumn{2}{|c|}{Preconditioner} & \multicolumn{3}{c|}{$\emptyset$} &  \multicolumn{3}{c|}{HypreBoomerAMG} & \multicolumn{3}{c|}{HypreILU}  \\
\hline  \multicolumn{2}{|c|}{CPU cores} & $256$ & $512$ & $1024$ & $256$ & $512$ & $1024$ & $256$ & $512$ & $1024$  \\
\hline  Total - Solver steps & Time (s) & $205.9$ & $242.7$ & $504.2$ & $264.4$ & $341.5$ & $\emptyset$ & $192.4$ & $199.2$ & $306.4$  \\
\hline  \multirow{2}*[1em]{\makecell[t]{Second to last\\ AMR iteration -\\ Last solver step}} & Number of iterations for PCG & $2580$ & $2571$ & $2352$ & $36$ & $36$ & $42$ & $890$ & $974$ & $1039$  \\
\cline{2-11}  & Time (s) & $16.0$ & $27.2$ & $58.1$ & $29.6$ & $37.9$ & $44.2$ & $13.8$ & $20.6$ & $22.4$  \\
\hline  \multirow{2}*[1em]{\makecell[t]{Last AMR\\ iteration -\\ Last solver step}} & Number of iterations for PCG & $4487$ & $3815$ & $3810$ & $38$ & $37$ & $\emptyset$ & $1038$ & $1059$ & $1253$  \\
\cline{2-11}  & Time (s) & $61.4$ & $49.3$ & $111.6$ & $57.6$ & $72.0$ & $\emptyset$ & $39.2$ & $29.8$ & $67.1$  \\
\hline 
\end{tabular}
\caption{Solving times - 3D Hertzian contact - AMR Combination 1
  ($e_{\Omega} = 1.5 \%$) and HyprePCG solver with different
  preconditioners - number of final mesh conforming DOFs: $9.0 \times 10^{6}$ - 6th-order transformation}
\label{tab:times_combinaison1}
\end{table}

As said before, the simulation with HypreBoomerAMG preconditioner
fails to converge with 1024 cores. This is due to a
non-symmetric-positive-definite (non SPD) preconditioning matrix
issue. This type of problem was also encountered by modifying the set
of parameters of the study case, for example by using AMR Combination
2 or by changing the value of the penalty coefficient. Modifying the
Strong Threshold, the HypreBoomerAMG parameter defined as the most
important option by MOOSE platform developers~\cite{INL-2023}, does
not enable us to tackle this issue. Since using the conjugate gradient
solver developed in MFEM leads to exactly the same issue, it seems
that the problem
comes from the association of this preconditioner with a conjugate
gradient solver, but not from the solver itself. In addition,
simulations with HypreILU as a preconditioner or without a
preconditioner are going well. During our studies, we observed that
using solvers that don't require a symmetric-positive-definite (SPD)
preconditioning matrix, such as GMRES or BiCGSTAB, allows you to get
round this HyperBoomerAMG issue, at the cost of larger computational
time. Thus, although using HypreBoomerAMG preconditioner can
drastically reduce the number of solver iterations, it is not
interesting (in terms of total computation time) to combine it with
the HyprePCG solver for our case studies.  The use of a HypreILU-type
preconditioner reduces the number of solver iterations by a factor of
2 to 3 compared to the case without preconditioner, and also appears
to be powerful in terms of computation time reduction. Consequently,
its combination with HyprePCG solver is efficient and should be
preferred for our contact problems.

A loss of performance is also observed in
Figure~\ref{fig:allKernelCompPourcentage_amrglob2} between 16 and 32
CPU cores. A conceivable explanation is that the intra-node bandwidth
seen by the CPU core is higher for 16 cores than for 32.\\

For AMR Combination 2, similar scalability analyses were performed and
are reported in Figure~\ref{fig:allKernelCompPourcentage_amrloc2}. The
threshold values $e_{\Omega, \text{LOC}} = 4\%$ and $\delta = 0.1\%$
have been set to highlight the specificity of this combination of
criteria for the same number of refinement steps as previously for AMR
Combination 1. The observations and conclusions regarding the global scalability
are the same as for the previous choice of AMR criteria. However,
the MFEM function forming the linear system performs better than
ideally for 128, 256 and 512 CPU cores, which is certainly due to the
presence of cache effects as we use much more cache memory with an
increasing number of nodes.
Indeed, if an increasing number of cores are used, then it is typically possible to employ more and more fast caches simultaneously.
In our case, on the topaze supercomputer, the L2 cache is 512 KB per core.
Thus, this means there is more room for data to fit in the L2 cache for an increasing number of cores.
In a strong scaling case, the total amount of data that can fit in the L2 cache increases with the number of cores.
The superlinear speedups reported can be explained by these cache effects. More specifically, the hypre library routine \texttt{hypre\_BoomerAMGBuildCoarseOperator} is responsible for this superlinearity, and consists in performing the matrix operation defined by Equation~\eqref{eqn:bilinear_form}.
In addition, let us underline that,
as the final number of conforming DOFs is almost half that for the previous AMR Combination 1 simulation, the total times are reduced accordingly.

\section{Conclusions}

In this paper, a scalable non-conforming AMR strategy dedicated to
elastostatic contact mechanics problems has been introduced. It has
been developed in the MFEM library. It is based on the implementation
of a node-to-node pairing combined with a penalization technique for
contact solution. An efficient partitioning that gathers paired
contact nodes together in the same parallel region is proposed. The
MFEM parallelized algorithm for nonconforming h-refinement has been
enriched by our own ESTIMATE-MARK-REFINE strategy. This strategy is in
particular designed to be effective for contact problems with curved
boundaries and quadrilateral/hexahedral first-order finite elements
thanks to the use of super-parametric finite elements. Representative
academic 2D and 3D examples confirm the accuracy of the proposed
approach.

The proposed strategy makes it possible to combine both efficient
parallelization and AMR techniques to solution 3D contact
problems. This is an important novelty as the combination of these
advanced numerical strategies was not previously considered. The
robust strategy we propose is scalable across hundreds of processes
and enables the treatment of problems with tens of millions of
unknowns. Ongoing studies will focus on the extension to contact
treatment of classical geometric partitioning algorithms
(Space-Filling Curve) and graph partitioning algorithms (using Metis
software package~\cite{Karypis-1998}). A enhanced gain in scalability
is expected with these partitionings.

Our strategy is aimed to be extended to more generic 3D contact
mechanics solutions, enabling a wider range of problems to be
considered. For example, an ongoing work is to implement
node-to-surface pairing. The main issue relies on the definition of a
relevant mesh partitioning. To solve industrial problems of interest,
our strategy has also to be extended to non-linear mechanics. The
MFEM-MGIS project~\cite{Helfer-2025, Helfer-MFEMMGIS} may be used to this end, since
the current standalone MFEM library only deals with elastic or
hyperelastic mechanical behaviors. Finally time-dependent problems
have to be considered, with a special attention to field transfer (for
discretization error control) and load balancing (for computational
time) issues between time steps.



\section*{Acknowledgements}
This work was developed within the framework of the MISTRAL joint research
laboratory between Aix-Marseille University, CNRS, Centrale Marseille, and CEA (Commissariat à
l’Énergie Atomique et aux Énergies Alternatives).
The authors would like to acknowledge the High Performance Computing
Center of CCRT for supporting this work by providing scientific
support and access to computing resources. We also acknowledge
financial support from the PICI2 project (CEA).

\bibliographystyle{elsarticle-num} 
\bibliography{biblio_article_Contact_AMR_HPC}


\end{document}

%% file: figures/patate_contact.eps_tex
\begingroup%
  \makeatletter%
  \providecommand\color[2][]{%
    \errmessage{(Inkscape) Color is used for the text in Inkscape, but the package 'color.sty' is not loaded}%
    \renewcommand\color[2][]{}%
  }%
  \providecommand\transparent[1]{%
    \errmessage{(Inkscape) Transparency is used (non-zero) for the text in Inkscape, but the package 'transparent.sty' is not loaded}%
    \renewcommand\transparent[1]{}%
  }%
  \providecommand\rotatebox[2]{#2}%
  \newcommand*\fsize{\dimexpr\f@size pt\relax}%
  \newcommand*\lineheight[1]{\fontsize{\fsize}{#1\fsize}\selectfont}%
  \ifx\svgwidth\undefined%
    \setlength{\unitlength}{134.01723353bp}%
    \ifx\svgscale\undefined%
      \relax%
    \else%
      \setlength{\unitlength}{\unitlength * \real{\svgscale}}%
    \fi%
  \else%
    \setlength{\unitlength}{\svgwidth}%
  \fi%
  \global\let\svgwidth\undefined%
  \global\let\svgscale\undefined%
  \makeatother%
  \begin{picture}(1,1.95630736)%
    \lineheight{1}%
    \setlength\tabcolsep{0pt}%
    \put(0,0){\includegraphics[width=\unitlength]{patate_contact.eps}}%
    \put(0.58,1.43){\color[rgb]{0,0,0}\makebox(0,0)[lt]{\begin{minipage}{1.22041209\unitlength}\raggedright $\Gamma_{C}^{1}$\end{minipage}}}%
    \put(0.58,0.85){\color[rgb]{0,0,0}\makebox(0,0)[lt]{\begin{minipage}{1.22041209\unitlength}\raggedright $\Gamma_{C}^{2}$\end{minipage}}}%
    \put(0.75,1.88){\color[rgb]{0,0,0}\makebox(0,0)[lt]{\begin{minipage}{1.22041209\unitlength}\raggedright $\Gamma_{BC}^{1}$\end{minipage}}}%
    \put(0.67,0.25){\color[rgb]{0,0,0}\makebox(0,0)[lt]{\begin{minipage}{1.22041209\unitlength}\raggedright $\Gamma_{BC}^{2}$\end{minipage}}}%
    \put(0.65,1.28){\color[rgb]{0,0,0}\makebox(0,0)[lt]{\begin{minipage}{1.22041209\unitlength}\raggedright ${\color{red} \textbf{n}^{1}}$\end{minipage}}}%
    \put(0.65,0.96){\color[rgb]{0,0,0}\makebox(0,0)[lt]{\begin{minipage}{1.22041209\unitlength}\raggedright ${\color{red} \textbf{n}^{2}}$\end{minipage}}}%
    \put(-0.05,1.12){\color[rgb]{0,0,0}\makebox(0,0)[lt]{\begin{minipage}{1.22041209\unitlength}\raggedright ${\color{OliveGreen} d}$\end{minipage}}}%
    \put(0.58,1.7){\color[rgb]{0,0,0}\makebox(0,0)[lt]{\begin{minipage}{1.22041209\unitlength}\raggedright $\Omega^{1}$\end{minipage}}}%
    \put(0.58,0.5){\color[rgb]{0,0,0}\makebox(0,0)[lt]{\begin{minipage}{1.22041209\unitlength}\raggedright $\Omega^{2}$\end{minipage}}}%
  \end{picture}%
\endgroup%

%% file: figures/hertz_contact.eps_tex
\begingroup%
  \makeatletter%
  \providecommand\color[2][]{%
    \errmessage{(Inkscape) Color is used for the text in Inkscape, but the package 'color.sty' is not loaded}%
    \renewcommand\color[2][]{}%
  }%
  \providecommand\transparent[1]{%
    \errmessage{(Inkscape) Transparency is used (non-zero) for the text in Inkscape, but the package 'transparent.sty' is not loaded}%
    \renewcommand\transparent[1]{}%
  }%
  \providecommand\rotatebox[2]{#2}%
  \newcommand*\fsize{\dimexpr\f@size pt\relax}%
  \newcommand*\lineheight[1]{\fontsize{\fsize}{#1\fsize}\selectfont}%
  \ifx\svgwidth\undefined%
    \setlength{\unitlength}{192.79584749bp}%
    \ifx\svgscale\undefined%
      \relax%
    \else%
      \setlength{\unitlength}{\unitlength * \real{\svgscale}}%
    \fi%
  \else%
    \setlength{\unitlength}{\svgwidth}%
  \fi%
  \global\let\svgwidth\undefined%
  \global\let\svgscale\undefined%
  \makeatother%
  \begin{picture}(1,1.20863437)%
    \lineheight{1}%
    \setlength\tabcolsep{0pt}%
    \put(0,0){\includegraphics[width=\unitlength]{hertz_contact.eps}}%
    \put(0.5,0.883788555){\color[rgb]{0,0,0}\makebox(0,0)[lt]{\begin{minipage}{1.22041209\unitlength}\raggedright $\Omega^{1}$\end{minipage}}}%
    \put(0.5,0.32485){\color[rgb]{0,0,0}\makebox(0,0)[lt]{\begin{minipage}{1.22041209\unitlength}\raggedright $\Omega^{2}$\end{minipage}}}%
    \put(0.55,-0.01){\color[rgb]{0,0,0}\makebox(0,0)[lt]{\begin{minipage}{1.22041209\unitlength}\raggedright $u_{D}$\end{minipage}}}%
    \put(0.5,1.2614){\color[rgb]{0,0,0}\makebox(0,0)[lt]{\begin{minipage}{1.22041209\unitlength}\raggedright $- u_{D}$\end{minipage}}}%
    \put(0.52,0.6597){\color[rgb]{0,0,0}\makebox(0,0)[lt]{\begin{minipage}{1.22041209\unitlength}\raggedright $O$\end{minipage}}}%
    \put(0.482,0.6257){\color[rgb]{0,0,0}\makebox(0,0)[lt]{\begin{minipage}{1.22041209\unitlength}\raggedright $+$\end{minipage}}}%
    \put(0.19,0.6357){\color[rgb]{0,0,0}\makebox(0,0)[lt]{\begin{minipage}{1.22041209\unitlength}\raggedright $\delta_{0}$\end{minipage}}}%
  \end{picture}%
\endgroup%